\documentclass[11pt]{article}

\usepackage{amssymb}
\usepackage{graphicx}
\usepackage{color}
\usepackage{amsmath,amsthm,amsfonts,epsfig,setspace}
\numberwithin{equation}{section}

\newtheorem{thm}{Theorem}[section]
\newtheorem{rmk}{Remark}[section]
\newtheorem{cor}{Corollary}[section]
\newtheorem{definition}{Definition}[section]
\newtheorem{lem}{Lemma}[section]

\newtheorem{prop}{Propsition}
\newtheorem{asump}{Assumption}

 \def\p{\partial}
\def \Vh0{\stackrel{\circ}{V}_h} \def\to{\rightarrow}

\newcommand{\q}{\quad}   
\def\l{\label}  \def\f{\frac}  

\def\m{\mbox}

\def\l|{\left|}
\def\r|{\right|}

\newcommand{\R}{\mathbb{R}}

\newcommand{\lc}
{\mathrel{\raise2pt\hbox{${\mathop<\limits_{\raise1pt\hbox
{\mbox{$\sim$}}}}$}}}

\newcommand{\gc}
{\mathrel{\raise2pt\hbox{${\mathop>\limits_{\raise1pt\hbox{\mbox{$\sim$}}}}$}}}

\newcommand{\ec}
{\mathrel{\raise2pt\hbox{${\mathop=\limits_{\raise1pt\hbox{\mbox{$\sim$}}}}$}}}

\def\be{\begin{equation}} \def\ee{\end{equation}}

\def\bea{\begin{eqnarray}}  \def\eea{\end{eqnarray}}

\def\beas{\begin{eqnarray*}} \def\eeas{\end{eqnarray*}}

\def\bn{\begin{enumerate}} \def\en{\end{enumerate}}

\title{Stability for the lens rigidity problem}
\date{}

\author{Gang Bao\thanks{School of Mathematical Sciences, Zhejiang University, Hangzhou 310027, China; The research was supported in part by a Key Project of the Major Research Plan of NSFC (No. 91130004), a NSFC A3 Project (No.11421110002), NSFC Tianyuan Projects (No. 11426235; No. 11526211, and a  NSFC Innovative Group Fund (No.11621101). Email: {\bf drbaogang@gmail.com}}
\and Hai Zhang\thanks{Department of Mathematics, The Hong Kong University of Science and Technology, Clear Water Bay, Kowloon, Hong Kong, China. The research was supported 
by Hong Kong RGC grant 26301016 and a startup grant from HKUST.
Email: {\bf haizhang@ust.hk}}
}

\begin{document}
\maketitle

\begin{abstract}
Let $g$ be a Riemannian metric for $\R^d$ ($d\geq 3$) which differs from the Euclidean metric only in a smooth and strictly convex bounded domain
$M$. The lens rigidity problem is concerned with recovering the metric $g$ inside $M$ from the corresponding lens relation on the boundary $\p M$. In this paper, the stability of the lens rigidity problem
is investigated for metrics which are a priori close to a given non-trapping metric satisfying ``strong fold-regular'' condition. A metric $g$ is called strong fold-regular if for each point $x\in M$, there exists a set of geodesics passing through $x$ whose conormal bundle covers $T^*_{x}M$. Moreover, these geodesics contain either no conjugate points or only fold conjugate points with a non-degeneracy condition. Examples of strong fold-regular metrics are constructed. Our main result
gives the first stability result for the lens rigidity problem in the case of anisotropic metrics which allow conjugate points. The approach is based on the study of the linearized inverse problem of recovering a metric from its induced geodesic flow, which is a weighted geodesic X-ray transform problem for symmetric 2-tensor fields. A key ingredient is to show that
the kernel of the X-ray transform on symmetric solenoidal 2-tensor fields is of finite dimension.
It remains open whether the kernel space is trivial or not.
\end{abstract}

\medskip


\section{Introduction and statement of the main result}

Let $(M, g)$ be a compact Riemannian manifold with boundary $\p M$. Let $\mathcal{H}^{t}(g)$ (or $\mathcal{H}^{t}$ for the ease of notation) be the geodesic flow on the tangent bundle $TM$ and let
$SM$ be the unit tangent bundle. Denote
\beas
S_{+} \p M&=& \{(x, \xi): x\in \p M, \xi\in T_xM, |\xi|_{g}=1, \langle \xi, \nu(x) \rangle > 0 \}; \\
S_{-} \p M&=& \{(x, \xi): x\in \p M,  \xi\in T_xM, |\xi|_{g}=1, \langle\xi, \nu(x)\rangle < 0 \},
\eeas
where $\nu(x)$ is the unit outward normal at $x$ and $\langle \cdot, \cdot \rangle$ stands for the inner product. We now define the
lens relation. Heuristically, the lens relation encodes all the information about the geodesic flow one can obtain from the boundary. It contains not only the lengths of the unit speed geodesics connecting boundary points but also the entrance and exit directions. 
In the simplest form, we assume that the metric $g$ is known on the boundary.
For each $(x, \xi) \in S_{-} \p M$, define $L(g)(x, \xi)>0$ to be the first positive moment at which the unit speed geodesic passing through $(x, \xi)$ hits the boundary $\p M$ ($L(g)(x, \xi)=\infty$ is allowed in which case we call the corresponding geodesic trapped), and 
$\Sigma(g):  S_{-} \p M \to \overline{S_{+} \p M}$ by
\be  \label{s-relation}
\Sigma(g)(x,  \xi) = \mathcal{H}^{L(g)}(x, \xi).
\ee
The map $\Sigma(g)$ is called the scattering relation. The lens relation induced by the metric $g$ is defined to be the pair $(\Sigma(g), L(g))$. The manifold $(M, g)$ is said to be non-trapping if there exists $T>0$
such that $L(g)(x, \xi)\leq T$ for all $(x, \xi)\in S_{-} \p M$.

Throughout the paper, we shall restrict ourselves to the case where the metrics are known on the boundary and use the above definition for the lens relation. For the more general definitions that remove the requirement of knowing $g$ on $\p M$, we refer to \cite{SU09}.


The lens rigidity problem is concerned with recovering the metric $g$
from its induced lens relation. It has important applications in geophysics where 
people are interested in finding the inner structures of earth from measurements of elastic waves on the surface. In a simplified model, the inner structure can be represented by wave speed (or a Riemaniann metric in geometric settings).  Usually, two type of measurement can be made, one is the wave field, and the other is the travel time (the time takes for wave to propagate from one place to another). We refer to \cite{U03, BZ12} for the relations between these two types of measurements. The former measurement yields the inverse problem of recovering wave speed from boundary dynamic Dirichlet-to-Neumann map (DDtN) and the latter the travel time tomography. It was shown in \cite{BZ12, thesis-zhang} that the boundary DDtN map
is very sensitive to wave speed, namely, small perturbations in wave speed can leads large deviations to their associated DDtN maps. As a consequence, wave speed inversion by DDtN map is inefficient in some sense. On the other hand, the inversion by the travel time has Lipshtiz type stability and hence can yield good and stable reconstruction. Acutally, the first successful application of travel time tomography in geophysics was obtained by Herglotz \cite{Herglotz05}, Wiechert and K.Zeoppritz in \cite{W07}, in which the wave speed inside the earth was assumed to be isotropic and spherically symmetry. For more general anisotropic wave speed, only travel time may not be sufficient to guarantee an accurate and stable reconstruction, and one has to use additional information from the measurement of wave fields. One natural information is the scattering relation. This is one of the motivations for studying the lens rigidity problem. 
We refer to \cite{S94, S99} and the references therein for the other developments in travel time tomography and the related problems in integral geometry. 
The lens rigidity problem is also closely related to the inverse scattering problem for metric perturbations of the Laplacian \cite{Gu76, E98} in the Euclidean space. In addition, it is also considered in the study of the 
AdS/CFT duality and holography \cite{PR04}.

It is known that there is no uniqueness to the lens rigidity problem in general. The first non-uniqueness example comes from the diffeomorphisms which leave the boundary $\p M$ fixed. More precisely, if $\phi: M\to M$ is a diffeomorphism of $M$ such that $\phi|_{\p M}=Id$ and the pull-back of the metric $g$ under $\phi$, denoted by $\phi^*g$, coincides with $g$ on $\p M$, then the two metrics $g$ and $\phi^*g$ have the same lens relation, i.e. $(\tilde{\Sigma}(\phi^*g), \tilde{L}(\phi^*g))=(\tilde{\Sigma}(g), \tilde{L}(g))$.
In addition, trapping of geodesics also prevents the uniqueness of the problem.
We refer the readers to \cite{CK94} for some counterexamples.
Therefore, a natural formulation for
the lens rigidity problem is as follows \cite{U03, SU09}:
given a compact non-trapping Riemannian manifold $(M, g)$ with boundary $\p M$,
is the metric $g$ uniquely determined by
its induced lens relation, up to the actions of diffeomorphisms which leave the boundary fixed? 

For the lens rigidity problem, only a few results are available.
Croke \cite{Croke05} showed that the finite quotient space of a lens rigid manifold is lens rigid.
Stefanov and Uhlmann \cite{SU09} proved uniqueness up to diffeomorphisms fixing the boundary
for metrics sufficiently close to a generic regular metrics. Vargo \cite{V09} showed that a class of analytic metrics are lens rigid. On the other hand, some interesting special cases of the lens rigidity problem with trapped geodesics are studied in \cite{Croke13, CH}.
More recently, Stefanov, Uhlmann and Vasy \cite{SUV13} have proved the lens rigidty uniqueness and stability in the conformal class where the manifold can be foliated by strictly convex hypersurfaces. Their approach is based on the study of the local boundary rigidity problem of determining the conformal factor of a Riemannian metric near a strictly convex boundary point from the distance function measured at pairs of points nearby.
We also remark that
a variant of the lens rigidity problem when the metrics are restricted to a conformal class has been investigated in \cite{thesis-zhang} with a Lipschitz stability result.

To our best knowledge, no stability result on the lens rigidity problem is available for anisotropic metrics which allow conjugate points. In this paper, we aim to address this issue and give the first result. We remark that a closely related problem for the lens rigidity problem is the boundary rigidity problem, which concerns with the unique determination of a simple metric from its induced boundary distance function, i.e. the lengths of geodesics connecting boundary points.  Recall that a compact Riemannian manifold with boundary is called simple if its boundary is strictly convex and the geodesics in it have no conjugate points.
It was observed by Michel \cite{Mi81} that the lens rigidity and boundary rigidity problem are equivalent for simple metrics.
We refer to \cite{ Muk77, Muk82, MuR78, BG80,  B83, G83, BCG95, Otal90, Croke90, Croke91, CDS00, S94, E98, SU98, LSU, PU05, SU05, BI10} and the references therein for various uniqueness and stability results for the boundary rigidity problem.

Compared to the boundary rigidity problem, the main difficulty of the lens rigidity problem lies in the presence of conjugate points on the geodesics which are excluded in the former case. In \cite{SUV13}, this difficulty was bypassed by the following approach: first, a local version of the problem in a small neighbourhood of a strictly convex boundary point was considered. Since the geodesics connecting two boundary points there have no conjugate points, the problem becomes equivalent to a local version of the boundary rigidity problem; second, a layer stripping argument is applied to yield a global result. In order for this method to work, a foliation condition is assumed.

In this paper, we overcome the difficult of conjugate points by another approach. We consider the linearized problem 
which is a geodesic X-ray transform for tensor fields, and analyze the effects of conjugate points on the transform. We impose conditions
on the conjugate points to ensure a
stability estimate for the X-ray transform. This estimate then yields the stability estimate for the nonlinear problem. 

Our approach is motivated by the results of Stefanov and Uhlmann \cite{SU04, SU05, SU09} on the boundary rigidity problem and the uniqueness of the lens rigidity problem. However, ours are different in the following two aspects: 
(1) The stability result in \cite{SU05} and the uniqueness results in \cite{SU04, SU09} are based on linearizing the operator which maps metrics to the lengths of geodesics connection boundary points, while our approach is based on linearizing the operator which maps metrics to their associated geodesic flows. There are two main advantages for the latter: it allows geodesics which have conjugate points; besides, the linearization is more straightforward, though it requires a larger manifold to work with. We remark that the information used from the geodesic flow is in fact equivalent to the lens relation (see Appendix A). (2) The geodesic X-ray transforms of the linearized operators in \cite{SU04, SU05, SU09} assume no conjugate points, while our transform allows the presence of conjugate points. In addition, our transform has a variable weight as opposed to the constant weight in theirs.
Note that fold type conjugate points occur commonly on geodesics in a general Riemannian manifold. In fact, it is shown that in the set of all geodesics passing through a given point $x\in M$, the set of geodesics with fold type conjugate points has the same dimension as the whole set, and the set of geodesics with other types of conjugate points has a lower dimension. Therefore, it is necessary to consider geodesics with fold type conjugate points in the study of the X-ray transform in a general Riemannian manifold \cite{BZ12}.
In order to control the effects of fold conjugate points,
we introduce the concept ``strong fold-regular'', which was motivated by Stefanov and Uhlmann's result in \cite{SUAPDE}. By allowing the usage of geodesics with strong fold-regular conjugate points in the inversion for the X-ray transform, we are able to study the X-ray transform on a quite general class of non-simple manifolds.

Throughout the paper, we restrict our discussion to $\R^d$ instead of general Rimannian manifolds
for the ease of exposition. We assume that $M$ is a strictly convex domain in $\R^d$, and that the metrics equipped to $M$ are
identical to the Euclidean metric, denoted by $e$, in a neighbourhood of $\p M$. This requirement saves us from technicalities caused from the boundary and allows us to focus on the essential difficulty of the problem, which is due to the presence of conjugate points on geodesics.

Our main result on the lens rigidity problem reads as follows (the proof will be given in Section \ref{sec-final}).

\begin{thm} \label{thm-lens}
Let $g$ be a smooth Riemannian metric equipped to $M$, which is a smooth and strictly convex bounded domain in $\R^d$ with $d\geq 3$.
Let $\hat{M}\Subset M$.
Assume that: (1). the support of $g-e$ is contained in $\hat{M}$; (2).
the geodesic flow induced by $g$ is non-trapping in $M$ over time $T$; (3). the metric $g$ is strong fold-regular (see Section \ref{sec-mainresult-x-ray} for its definition).
Then there exist a positive integer $k$, a finite dimensional space $\mathfrak{L} \in L^2(S(\tau_2M))$ where $S(\tau_2M)$ denotes the set of symmetric 2-tensor fields on $M$, and a finite number of smooth functions $\alpha_j \in C^{\infty}(S_{-}\p M)$, $j=1, 2, ...N$,
such that for any Riemannian metric $\tilde{g}$ close enough to $g$ in $C^{k, \f{1}{2}}(M)$ and differs from $e$ only  in $\hat{M}$, there exists a
diffeomorphism $\varphi$ of $M$ with the property that
$\varphi^*\tilde{g}= e$ in a neighborhood of $\p M$, and moreover
\be \label{estimate00}
\|\varphi^*\tilde{g}- g\|_{C^{0}(M)} \lesssim  \sum_{j=1}^N \|\alpha_j (\mathcal{H}^T(\tilde{g}) -\mathcal{H}^T(g)) \|_{H^{3+[\f{d}{2}]}(S_{-}\p M)},
\ee
provided $\varphi^*\tilde{g}- g \perp \mathfrak{L}$ in $L^2(S(\tau_2M))$.
 \end{thm}

\begin{rmk}
It is shown in the proof of Theorem \ref{thm-lens} that the integer $k$ can be chosen to be $13+2[\f{d}{2}]$. 

\end{rmk}

\begin{rmk}
In the right hand side of the estimate (\ref{estimate00}), the difference is between two vectors in $\R^{2d}$. This makes sense because of our assumption that the metrics considered here are Euclidean outside $M \subset \R^d$. For a general Riemannian manifold, we need to use the original lens relation (\ref{s-relation}) which is defined on $\p M$. The difference can be realized by using local coordinates near $\p M$, such as in \cite{SUV13}. 
\end{rmk}

\begin{rmk}
In Theorem \ref{thm-lens}, only partial information from the geodesic flow (or the lens relation) is used. As a consequence, the assumptions on the background metric $g$ can be significantly weakened in the following two ways:

(1). The non-trapping condition can be replaced by the condition that only
the geodesics which are used for the stability estimates (or those which are in the support of the cut-off functions $\alpha_j$'s) are non-trapping in $M$ for time $T$.

(2). The convexity condition on the boundary $\p M$ can be replaced by the condition that the geodesics used are transverse to the boundary $\p M$.
\end{rmk}

\begin{rmk}
In Theorem \ref{thm-lens}, if we have another metric $\tilde{g}_1$ which satisfies the same conditions as $\tilde{g}$, then the conclusion becomes
\[
\|\varphi^*\tilde{g}- \varphi^*_1\tilde{g}_1\|_{C^{0}(M)} \lesssim  \sum_{j=1}^N \|\alpha_j (\mathcal{H}^T(\tilde{g}) -\mathcal{H}^T(\tilde{g}_1)) \|_{H^{3+[\f{d}{2}]}(S_{-}\p M)},
\]
provided $\varphi^*\tilde{g}- \varphi^*_1\tilde{g}_1 \perp \mathfrak{L}$, where $\varphi_1$ is determined by $\tilde{g}_1$ in the same way as $\varphi$ by $\tilde{g}$.
\end{rmk}

\begin{rmk}
The finite dimensional space $\mathfrak{L}$ is related to the kernel of the linearized inverse problem of determining $g$ from
its induced geodesic flow $\mathcal{H}^T(g)|_{S_{-}\p M}$, or more precisely, the geodesic X-ray transform $\mathfrak{X}$ (see Section \ref{sec-linearization}) for solenoidal 2-tensor fields.
It is not clear that whether $\mathfrak{L}$ is trivial or not for a general non-simple metric. We note that the methods developed in this paper is unlikely to yield an answer to this question. 
For simple metrics, the geodesic X-ray transform is defined by linearizing the boundary distance function, and the corresponding problem is referred to as the ``s-injectivity problem'', see for instance \cite{SU04, SU08AJM, GMG13}. 
It was proven that s-injectivity held for metrics with negative curvature in \cite{PS88},  for 
metrics with small curvature and in \cite {S94},
for Riemannian surfaces with no focal points in \cite{ShU01}, and for generic analytical metrics
in \cite{SU05}. 
\end{rmk}

\begin{rmk}
Theorem \ref{thm-lens} can be stated in terms of the lens relation instead of the geodesic flow. In fact, there are two ways to do so. One way is to replace the geodesic flow by the corresponding lens relation with the help of Identity (\ref{identity3-1-1}). The other way is to linearize the lens relation with respect to the metric as for the geodesic flow. We can show that a similar conclusion holds if we replace $\mathcal{H}^T$ by $(\Sigma, L)$.
\end{rmk}


We now give a brief account of the main approach in the paper.
We first linearize the operator which associates a metric with its induced geodesic flow on the cosphere bundle at a fixed smooth background metric $g$. We assume that $g$ is non-trapping in $M$ over the time $T$. This leads to an X-ray transform operator $\mathfrak{X}$ for symmetric 2-tensor fields.
We form the normal operator $\mathfrak{N}=\mathfrak{X}^{\dag}\mathfrak{X}$. In the case when the metric $g$ is not simple (or the exponential maps of the manifold $(M, g)$ have singularities), the Schwartz kernel of $\mathfrak{N}$ has two types of singularities: one is from the diagonal, the other is from conjugate points. By the arguments in \cite{SUAPDE}, we show that the former singularities yield a $\Psi$DO (pseudo differential operator) while the latter ones yield a FIO (Fourier integral operator) locally in the case when the singularities are of fold type.

We study the stability of the above X-ray transform operator $\mathfrak{X}$ for symmetric 2-tensor fields. Since it vanishes for potential 2-tensor fields, we restrict to the solenoidal 2-tensor fields. By micro-local analysis, a set of geodesics passing through $x$ whose conormal bundle can cover the cotangent space $T_x^*\R^d$ are needed in order to reconstruct the solenoidal part of a 2-tensor field $f$ at $x$ from its X-ray transform $\mathfrak{X}f$. We allow fold conjugate points along these geodesics, but require that these conjugate points contribute to a smoother term than the point $x$ itself. This is the case when the intrinsic derivatives of the differential of the exponential map at the conjugate vectors satisfy a non-degenerate condition, which is called ``strong fold-regular''. We establish a local stability estimate near ``strong fold-regular'' points. These local estimates are glued together to form a global one. We remark that it remains open whether the X-ray transform operator has a non-trivial kernel in the space of solenoidal 2-tensor fields. This problem is often referred to as the ``s-injectivity'' problem.
By imposing some orthogonality conditions,  we obtain a Lipschitz type stability estimate (see Theorem \ref{thm-3}).

We remark that a similar X-ray transform operator is obtained by linearizing the operator which maps metrics within a conformal class to their induced geodesic flows at a smooth background metric in \cite{BZ12}.
The transform there is applied to scalar functions. Its kernel can be shown to be of finite dimension. In comparison, the X-ray transform $\mathfrak{X}$ in this paper has an infinitely dimensional kernel space. In fact, it vanishes for all potential 2-tensor fields as shown in Lemma \ref{lem-ns}. Thus the analysis for its invertibility and stability is much more subtle here. See also the previous works on this issue \cite{S99, CDS00, SU04, SU05, SU09}.

The Lipschitz type estimate for the X-ray transform $\mathfrak{X}$, or the linearized inverse problem of recovering a metric from its induced geodesic flow,
enables us to derive a Lipschitz stability result for the nonlinear inverse problem. This last step is made possible by using the construction for diffeomorphisms in \cite{SU09}.


The paper is organized as follows: In Section 2, we present some preliminaries and introduce two X-ray transform operators $\mathfrak{I}$ and $\mathfrak{X}$ which are obtained from linearizing the operator which maps a metric to its induced geodesic flow. The results on the stability estimates for the X-ray transform $\mathfrak{X}$ are stated in Section \ref{sec-mainresult-x-ray}. In Section \ref{sec-i-m}, we derive properties of the X-ray transform $\mathfrak{I}$ and its normal $\mathfrak{M}$. Based on these properties, we study the X-ray transform $\mathfrak{X}$ in Section \ref{sec-x-n} and prove the stability results stated in Section \ref{sec-mainresult-x-ray}. 
in Section \ref{sec-final}, we prove the main result of the paper, Theorem \ref{thm-lens}.
Finally, examples of strong fold-regular Riemannian metrics are constructed in Section \ref{sec-example}.

\section{Preliminaries}

\subsection{Notations}
Throughout the paper, we denote by $\exp_{x}(\cdot): T_x\R^d \to \R^d$ the exponential map at $x\in \R^d$ with respect to the background metric $g$. The differential of $\exp_{x}(\cdot)$ at
$\xi_*\in T_x\R^d$ is denoted by $d_{\xi} \exp_x(\xi_*)$. The norm of $\xi\in T_x\R^d$ is denoted by $|\xi|$. If the linear map $d_{\xi} \exp_x(\xi_*)$
is not a diffeomorphism, we call the vector $\xi_*$ a conjugate vector, and the corresponding point $\exp_x(\xi_*)$ a conjugate point. The kernel of the linear map $d_{\xi} \exp_x(\xi_*)$ is denoted by $N_x(\xi_*)$, which is a subspace of $T_x\R^d$.

\medskip

In addition, we introduce the following conventions:
\bn

\item
Let $A_1$ and $A_2$ be two matrices (including vectors which can be regarded as single column or single row matrices), then the product of $A_1$ and $A_2$ is denoted by $A_1\cdot A_2$. Sometimes, the dot is omitted for simplicity;

\item
Let $A$ be a complex matrix, then $A^{\dag}$ stands for its conjugate transpose. If $A$ is a linear operator in a Hilbert space, then $A^{\dag}$ stands for its formal adjoint. If $A$ is real and symmetric and $C$ a real number, then $A\geq C$ means that the
matrix $A-C\cdot Id$ is symmetric and positive definite.

\item
Let $U$ and $V$ be two open set in a metric space, then $U\Subset V$ means that the closure of $U$, denoted by $\bar{U}$ is compact and is a subset of $V$;

\item
Let $C_1$ and $C_2$ be two positive numbers, then $C_1 \lesssim C_2$ means that $C_1\leq C\cdot C_2$ for some constant $C>0$ independent of $C_1$ and $C_2$.
\en

\subsection{Symmetric tensor fields}
We denote by $S(\tau_2M)$ the set of symmetric covariant 2-tensor fields to the manifold $M$. In the natural coordinate of $\R^d$, each $f\in S(\tau_2M)$ is assigned a
family of functions $\{f_{ij};  1\leq i, j \leq d \}$ in $M$ such that $f_{ij}=f_{ji}$.
We denote by $S(\tau_2M, \R^d)$ the set of families of functions $\{\Pi_{ij}^k; 1\leq i, j, k \leq d \}$ in $M$ such that $\Pi_{ij}^k=\Pi_{ji}^k$.
An example of $S(\tau_2M, \R^d)$ is the set of Christoffel symbols. We remark that elements in $S(\tau_2M, \R^d)$ are not tensor fields, i.e. they are not invariant under the change of coordinates. We also remark that we shall use the natural coordinate of $\R^d$ throughout the paper.

Let $L^2(S(\tau_2M))$ be the Hilbert space of $L^2$-integrable symmetric covariant 2-tensor fields with the following inner product
\be  \label{inner-product}
(f, h)_{L^2(M)}= \int_{M} f_{ij} h^{ij} \sqrt{\det g} dx = \int_{M} f_{ij} g^{ii'}g^{jj'}h_{i'j'} \sqrt{\det g} dx,
\ee
where $\{h^{ij}\}$ is the contravariant tensor field corresponding to $\{h_{ij}\}$.

Similar to the above $L^2$ space for 2-tensor fields, we can define $L^2$ space for any tensor fields.
Here we remark that throughout the paper,
the Einstein's summation rule is applied whenever there are sub and superscripts with the same label. For convenience, we always use the usual convection of raising and lowering indices and we treat the covariant and the corresponding contravariant tensor field as two representations of the same tensor field.

Based on the $L^2$ space defined above, one can define the $H^k$ space for any positive integer $k$ by the following inductive formula:
\[
(f, h)_{H^k(M)} = (f, h)_{L^2(M)} + (\nabla f, \nabla h)_{H^{k-1}(M)}.
\]
By the standard interpolation theory, one can define the $H^s$ space for any $s\geq 0$ as well.

\begin{rmk}
One can define an equivalent norm for the space of 2-tensor fields $L^2(S(\tau_2M))$ by using the natural coordinate of $\R^d$. Indeed, for $f, h \in L^2(S(\tau_2M))$, define
\be  \label{inner-product2}
(f, h)_{\tilde{L}^2(M)} = \sum_{ij} \int_{M}f_{ij}(x)h_{ij}(x)dx.
\ee

We can check the equivalence of the two norms by using the open mapping theorem.
Similarly, we can also define the equivalent norm for the space $H^k(S(\tau_2M))$.

The norm induce by the inner product (\ref{inner-product2}) shall be
used frequently in the subsequent analysis, since it is more convenient to work with than the one induced by (\ref{inner-product}). Considering that the two norms are equivalent, we use the
same symbol $\|\cdot\|_{L^2}$ (or $\|\cdot\|_{H^k}$) to denote them for simplicity.
\end{rmk}

We define $L^2(S(\tau_2M), \R^d)$ to be the Hilbert space of $L^2$-integrable fields in $S(\tau_2M, \R^d)$ with the following inner product
$$
(\Pi, \Gamma)= \int_{M} \Pi_{ij}^k \Gamma_{k}^{ij} \sqrt{\det g} dx,
$$
where $\Gamma_{k}^{ij} = \delta_{kk'}g^{ii'}g^{jj'}\Gamma_{i'j'}^{k'}$ with $\delta_{kk'}$ being the standard Cronecker's symbol.
We define the $H^k$ space for $S(\tau_2M, \R^d)$
by using the inner product (\ref{inner-product2}) for $k\geq 1$.

Let $\tilde{M}$ be a compact set in $M$ and let $k$ be a nonnegative integer, we denote by $H^k(S(\tau_2(\tilde{M}, M)))$ the subspace of
$H^k(S(\tau_2M))$ whose elements are supported in $\tilde{M}$.

Finally, we introduce a lifting operator. For each $h \in S(\tau_2M, \R^{d})$, we define $\iota (h) \in S(\tau_2M, \R^{2d})$ by setting the first $d$ components to be zero.

\subsection{Decomposition of symmetric 2-tensor fields}
We present some basic facts about the decomposition of symmetric 2-tensor fields. We refer to \cite{S99, SU05} for details. Given a symmetric 2-tenor field $f=\{f_{ij}\}$, we define the 1-tensor field
$\delta^s f$ by
$$
 [\delta^s f]_i = g^{jk} \nabla_kf_{ij},
$$
where $\nabla_k$ is the covariant derivative. On the other hand, given a 1-tensor field $v=\{v_i\}$, we define the 2-tensor field $d^s v$, called symmetric differential of $v$, by
$$
 [d^sv]_{ij}= \f{1}{2}(\nabla_i v_j + \nabla_j v_i).
$$
The operators $d^s$ and $-\delta^s$ are formally adjoint to each other in the $L^2$ space.

For each $f\in L^2(S(\tau_2M))$, there exists a unique orthogonal decomposition
$$
f= f^s + d^s v = \mathcal{S}f + \mathcal{P}f,
$$
such that $f^s \in L^2(S(\tau_2M))$ satisfies $\delta^s f^s =0$ and $v\in H^1_0(M)$ (i.e. $v=0$ on $\p M$). The fields $f^s$ and $d^sv$ are called the solenoidal and potential parts of $f$, respectively. The following orthogonality property holds for
the two projectors $\mathcal{S}, \mathcal{P}:L^2(S(\tau_2M)) \to L^2(S(\tau_2M)) $:
$$
\mathcal{P} \mathcal{S} = \mathcal{S}\mathcal{P}=0.
$$

The decomposition can be constructed as follows. Consider the operator $\triangle^s_{M} = \delta^s d^s$. It can be verified that
$\triangle^s_M$ is elliptic in $M$ and the corresponding Dirichlet problem satisfies the Lopatinskii condition. Thus $\triangle^s_M$ is invertible with the Dirichlet boundary condition.  Denote by $(\triangle^s_{M})^{-1}$ its inverse. Then
$$
 v= (\triangle^s_{M})^{-1} \delta^s f, \quad f^s = f- d^s(\triangle^s_{M})^{-1} \delta^s f.
$$

%

\subsection{Decomposition of symmetric 2-tensors}
We introduce a similar decomposition for symmetric tensors (see Section 2.6 in \cite{S94} for detail). Denote
$$
S_2=\{\{v_{ij}\}; v_{ij}=v_{ji}\in \R, 1\leq i, j\leq d\}; \quad S_{1}=\R^d.
$$
For each $x\in \R^d \backslash 0$, we define $i_x: S_1 \to S_{2}$ and $j_x: S_2 \to S_{1}$ in the following way:
$$
(i_xv)_{ij}= \f{1}{2}(v_ix_j+ v_jx_i), \quad (j_xu)_i = u_{ij}x^j,
$$
where $v\in S_1$ and $u\in S_2$.
The two operators $i_x$ and $j_x$ are called the symmetric multiplication and convolution with $x$, respectively. We can show that for every $f\in S_2$ and  $x\in \R^d \backslash 0$, there exist uniquely determined $h\in S_2$ and $v\in S_{1}$ such that
$$
f= h+ i_xv ,\quad j_xh=0.
$$

\subsection{Linearization of the geodesic flow with respect to the metric } \label{sec-linearization}

We linearize the operator which maps a metric to its induced geodesic flow restricted to the unit sphere tangent bundle. A similar linearization in the cotangent bundle is done in \cite{SU98} and is used in the travel time tomography problem \cite{CQUZ11}.   
Recall that $g$ is the fixed background metric. Denote by $\Gamma= \{\Gamma_{ij}^k; 1\leq i, j,k \leq d\}$ its Christoffel symbol. Let $f\in C^{\infty}(S(\tau_2M))$ be such that
its support is contained in $\hat{M} \Subset M$. Then
$g+f$ defines a perturbed metric in $M$ provided $\|f\|_{C^2(M)} \ll 1$.
Denote by $\mathcal{H}^T(g+f)$ the corresponding geodesic flow.
We now derive an explicit formula for
the linearized map $\f{\delta \, \mathcal{H}^T(g)}{\delta \, g}$ at $g$ below.

We first define a matrix $\Phi(x, \xi)$
for each $(x, \xi)$ that lies on the orbit of the set $S_{-}\p M$ under the geodesic flow $\mathcal{H}^t(g)$ ($t\geq 0$), see also \cite{BZ12}. Let
$(x_0, \xi_0)=\mathcal{H}^{L^{-}(g)(x, \xi)}(g)(x, \xi) \in S_{-}\p M$,
where $L^{-}(g)(x, \xi)$ is the first negative moment the geodesic orbit
$\mathcal{H}^t(g)(x, \xi)$
hits the boundary $S_{-}\p M$. Let
$\phi(t, x_0, \xi_0)$ be the solution of the following ODE system
$$
\dot{\phi}(t, x_0, \xi_0)= -\phi(t, x_0, \xi_0)A(\mathcal{H}^t(g)(x_0, \xi_0)), \q \phi(0, x_0, \xi_0)= Id,
$$
where $A(x, \xi) =(\f{\p H}{\p x}, \f{\p H}{\p \xi})$ is a $2d\times 2d$ matrix and $H$ is given by
$$
H(x, \xi)= (\xi^1, \xi^2,...,\xi^d , -\Gamma_{ij}^1\xi^i\xi^j, -\Gamma_{ij}^2\xi^i\xi^j,...,-\Gamma_{ij}^d\xi^i\xi^j)^\dag.
$$
The matrix $\Phi(x, \xi)$ is defined by
$$
\Phi(x, \xi)= \phi(T, x_0, \xi_0)^{-1}\phi(-L^{-}(g)(x_0, \xi_0), x_0, \xi_0).
$$

It is clear that $\Phi(\cdot, \cdot)$ is well-defined and is smooth in a sufficiently small neighborhood of any unit speed geodesic orbit if only the corresponding geodesic
is transverse to the boundary $\p M$.


With the matrix-valued weight $\Phi$ defined above, we define the X-ray transform operator
$\mathfrak{I}: C^{\infty}(S(\tau_2M); \R^{2d}) \to \mathcal{D}'(S_{-}\p M; \R^{2d})$ by
\bea
(\mathfrak{I} \Pi)(x_0, \xi_0)
&=& \int_0^T\Phi(x(s), \xi(s))\Pi_{ij}(x(s))\xi^i(s)\xi^j(s) \,ds  \nonumber\\
&=& \int_0^{L(g)(x_0, \xi_0)}\Phi(x(s), \xi(s))\Pi_{ij}(x(s))\xi^i(t)\xi^j(s) \,ds \label{formula-I}
\eea
where $(x(s), \xi(s))= \mathcal{H}^{s}(g) (x_0, \xi_0)$ and $\Pi_{ij}\xi^i\xi^j$ is viewed as a vector in $\R^{2d}$ with the k-th component given by $\Pi^k_{ij}\xi^i\xi^j$.

We 
introduce one more operator.
For each $f\in S(\tau_2M)$, we define $\mathcal{L}f \in S(\tau_2M, \R^d)$
$$
(\mathcal{L}f)_{ij}^k= \f{1}{2}g^{lk}\left(\f{\p f_{jl}}{\p x^i} +
\f{\p f_{il}}{\p x^j}- \f{\p f_{ij}}{\p x^l}\right) - g^{lk} \Gamma^{m}_{ij}f_{ml}.
$$
It is clear that $\mathcal{L}$ is a first order partial differential operator. We denote by $\mathcal{L}^{\dag}$ is adjoint.

Recall that $\iota$ is the shift operator defined at the end of Section 2.2. Denote $\mathfrak{X}= \mathfrak{I}\circ \iota \circ\mathcal{L}$. Then the following result shows  that
$\f{\delta \, \mathcal{H}^T(g)}{\delta \, g}|_{S_{-}\p M}= \mathfrak{X}$. 

\begin{prop}  \label{prop1}
The following estimate holds
\[
 \|\mathcal{H}^T(g+f) - \mathcal{H}^T(g) - \mathfrak{X}f\|_{C^{k}(S_{-}\p M)} \lesssim \|f\|_{C^{k+2}(M)}^2
\]
 for any nonnegative integer $k$.
\end{prop}

{\bf{Proof}}. See Appendix B.

\begin{rmk}
Formula (\ref{formula-I}) is derived in the coordinate of $T\R^d$. It is not geometrically invariant.
\end{rmk}

\section{Stability for the geodesic X-ray transform $\mathfrak{X}$} \label{sec-mainresult-x-ray}

 ``Fold-regular'' metrics was introduced in \cite{BZ12} to study the geodesic X-ray transform for scalar functions in a general Riemannian manifold. They are motivated by the work \cite{SUAPDE}, and generalize the ``regular'' metrics in \cite{SU08AJM} by allowing fold conjugate points along the geodesics which are used in the inversion of the X-ray transform. In this section, we present some stability results for the transform $\mathfrak{X}$ obtained in the previous section
under the ``strong fold-regular'' condition. We refer to \cite{S94, S99, SU04, SU05, GMG13} and the references therein for the study of X-ray transforms in simple manifolds and their applications. We also refer to \cite{GV} for the study of X-ray transform for scalar functions in a class of nonsimple manifolds which satisfy a global foliation condition.

We first introduce some definitions.
We refer to \cite{AZV, SUAPDE, BZ12} for more more details and discussions.
Recall that a conjugate vector $\xi_*\in T_x\R^d$ is of fold type if the following two conditions are satisfied: (1) the rank of the linear mapping $d_{\xi} \exp_x(\xi_*)$ equals to $d-1$ and the function $\det (d_{\xi}\exp_x(\xi))$ vanishes of order 1 at $\xi_*$; (2) the kernel space $N_x(\xi_*)$ of the linear mapping $d_{\xi} \exp_x(\xi_*)$ is transversal to the manifold $\{\eta: \det (d_{\eta} \exp_x(\eta))=0\}$ at $\xi_*$.

\begin{definition}
A fold vector $\xi_*\in T_{x}M$ is called strong fold-regular if the following condition is satisfied
\be \label{graph-condition}
  d^2_{\xi} \exp_{x} (\xi_*)(N_{x}(\xi_*)\setminus 0 \times \cdot )|_{T_{\xi_*}S(x)} \quad \mbox{is of full rank}
\ee
where
$S(x)= \{ \xi \in T_{x} \R^d; \, \det (d_{\xi} \exp_{x} \xi)=0\}$ and $T_{\xi_*}S(x)$ is its tangent space at $\xi_*$.
\end{definition}

The above condition (\ref{graph-condition}) was first introduced in \cite{SUAPDE}, which
guaranties the graph condition (\cite{H-III85}) for the FIO obtained from the X-ray transform studied therein.
We use it for the same purpose here, see Lemma \ref{lem-modified-n2}.

\begin{definition} \label{def-fold-regular-point}
A point $x\in M$ is called strong fold-regular if there exists a compact subset $\mathcal{Z}(x) \subset S_{x}M$ such that the following two conditions are satisfied:
\begin{enumerate}
\item
For each $\xi \in \mathcal{Z}(x)$, there exist either no singular vectors or those of strong fold-regular type along the ray $\{t\xi:\, t\in \R\}$ for the map
$\exp_x (\cdot)$ before it hits the boundary; moreover, the corresponding geodesic hits the boundary $\p M$ transversely.

\item
$\forall \, \xi \in S_{x}M, \, \exists \, \theta \in \mathcal{Z}(x) , \mbox{ such that }\, \theta \perp \xi$.
\end{enumerate}
\end{definition}

The second condition in the above definition can be viewed as ``completeness'' condition, see \cite{SU08AJM} for instance. It ensures the ellipticity of the $\Psi$DO part of the normal of properly truncated X-ray transforms,
see Lemma \ref{lem-n111}.


We next introduce the truncated X-ray transform operator.
For any $\alpha\in C^{\infty}_0(S_{-}\p M)$, we define the truncated operator $\mathfrak{X}_{\alpha}$ by
\begin{equation} \label{formula-X-alpha}
\mathfrak{X}_{\alpha}f(x_0, \xi_0)
= \alpha(x_0, \xi_0) \mathfrak{X}f(x_0, \xi_0).
\end{equation}

Denote by $\mathfrak{N}_{\alpha}= \mathfrak{X}_{\alpha}^{\dag}\mathfrak{X}_{\alpha}$ the normal operator.

We now present a local stability estimate near strong fold-regular points and several corollaries.  The proofs will be given in Section \ref{sub-sec-proof-x-ray}.

\begin{thm} \label{thm-1}
Let $x_*\in \tilde{M}\Subset M$ be a strong fold-regular point and let $k$ be a nonnegative integer. Then there exist two neighborhoods
$U(x_*)\Subset \tilde{U}(x_*)$ of $x_*$, and a cutoff function
$\alpha \in C_0^{\infty}(S_{-}\p  M)$, such that the following estimate holds uniformly for all
$f\in H^k(S(\tau_2(\tilde{M}, M)))$
\be \label{estimate1}
 \|f^s\|_{H^k(U(x_*))} \lesssim \|\mathfrak{N}_{\alpha}f\|_{H^{k-1}(\tilde{U}(x_*))} +  \|f\|_{H^{t}({M})},
\ee
where $t= \max \, \{k+1-\f{d}{2}, k-1\}$.
\end{thm}

The above local result can be extended to a global one.

\begin{definition}
The background metric $g$ is called strong fold-regular if all points in $M$ are strong fold-regular with respect to the geodesic flow $\mathcal{H}^t(g)$.
\end{definition}

\begin{cor} \label{thm-2}
Assume that the background metric $g$ is strong fold-regular. Let $\tilde{M}$ be a strictly convex and smooth sub-domain of $M$. Let $k$ be a nonnegative integer. Then there exist
$\tilde{U}(x_j) \subset M$, $\alpha_j \in C_0^{\infty}(S_{-}\p  M)$, $j=1, 2,...,N$, such that
the following estimate holds for all $f\in H^k(S(\tau_2(\tilde{M}, M)))$:
\[
 \|f^s\|_{H^k(M)} \lesssim \sum_{j=1}^N \|\mathfrak{N}_{\alpha_j}f\|_{H^{k-1}(\tilde{U}(x_j))} +  \|f\|_{H^{t}({M})}.
 \]
 where $t= \max \, \{k+1-\f{d}{2}, k-1\}$
\end{cor}


\begin{cor} \label{thm-3}
Assume that the background metric $g$ is strong fold-regular. Let $\tilde{M}$ and $\alpha_j \in C_0^{\infty}(S_{-}\p  M)$, $j=1, 2,...,N$, be as in corollary \ref{thm-2}. Then there exists a
finite dimensional space $\mathfrak{L} \in \mathcal{S}L^2(S(\tau_2(\tilde{M}, M))) $
such that the following estimate holds for all $f\in L^2(S(\tau_2(\tilde{M}, M)))$ satisfying $f^s\perp \mathfrak{L}$ in $L^2(S(\tau_2M ))$:
\be  \label{equation41}
 \|f^s\|_{L^2(M)} \lesssim \sum_{j=1}^N \|\mathfrak{N}_{\alpha_j}f\|_{H^{-1}(\tilde{U}(x_j))}.
 \ee
\end{cor}

\begin{cor} \label{thm-4}
Assume that the background metric $g$ is strong fold-regular. Let $\tilde{M}$, $\mathfrak{L}$ and $\alpha_j \in C_0^{\infty}(S_{-}\p  M)$, $j=1, 2,...,N$, be as in Corollary \ref{thm-3}, and let $k$ be an nonnegative integer. Then the following estimate holds for all $f\in H^k(S(\tau_2(\tilde{M}, M)))$ satisfying $f^s\perp \mathfrak{L}$ in $L^2(S(\tau_2M ))$:
\be  \label{equation42}
 \|f^s\|_{H^k(M)} \lesssim \sum_{j=1}^N \|\mathfrak{N}_{\alpha_j}f\|_{H^{k-1}(\tilde{U}(x_j))}.
 \ee
\end{cor}


\begin{rmk}
 The strong fold-regular condition in Corollary \ref{thm-3} and \ref{thm-4} may be weakened. In fact, as in \cite{BZ12},
we can define a fold vector $\xi \in T_xM$ to be fold-regular if the operator germ $\mathfrak{M}_{\xi}$, which characterizes contributions from an infinitesimal small neighborhood of the point $\exp_x\xi$ to the normal operator $\mathfrak{M}$, is compact from $H^k(\tilde{M}, M)$ to $H^{k+1}(U(x))$ for any integer $k$ and some neighborhood $U(x)$ of $x$. Then the results in Corollary \ref{thm-3} and \ref{thm-4} and consequently Theorem \ref{thm-lens} also hold under the fold-regular condition.

\end{rmk}

\section{The X-ray transform operator $\mathfrak{I}$ and its normal $\mathfrak{M}$} \label{sec-i-m}

We study the X-ray transform operator $\mathfrak{I}$ and its normal $\mathfrak{M}$ in this section.
We first introduce some measures. 
Let $x\in M$, we denote by $d\mu_x(\xi)$ and $d\mu(x,\xi)$ the measure induced by the metric $g$ to the sphere $S_xM$ and the spherical bundle $SM$, respectively. For the set $S_{-} \p M$ characterizing the set of geodesics passing through $M$, its canonical measure is given by $|\langle \nu(x), \xi \rangle|d\Sigma^{2d-2}(x, \xi)$, where $d\Sigma^{2d-2}(x, \xi)$ is the restriction of the measure $d\mu(x,\xi)$ to the subset $S_{-}\p M$ and $\nu(x)$ is the outward normal to $\p M$ at $x$.
In the coordinates of $\R^d$, we have
\beas
d\mu_x(\xi) &=& (\det{g})^{\f{1}{2}} \sum_{i=1}^d (-1)^{i-1} \xi^i d\xi^1\wedge...\wedge d\xi^{i-1}\wedge d\xi^{i+1}\wedge...\wedge d\xi^d; \\
d\mu(x,\xi) &=& d\mu_x(\xi) \wedge (\det{g})^{\f{1}{2}}dx.
\eeas
We also use $d\mu_x(\xi)$ to denote the measure to the tangent space $T_xM$. In coordinates,
$d\mu_x(\xi)=(\det{g})^{\f{1}{2}}d\xi $. We refer to \cite{S99} for more detail.

We now present some basic properties about the X-ray transform $\mathfrak{I}$.
\begin{lem}
 The X-ray transform operator $\mathfrak{I}$ is bounded from $H^k(S(\tau_2M, \R^{2d}))$ to $ H^k(S_{-}\p M, \R^{2d})$ for any integer $k\geq 0$.
\end{lem}

{\bf{Proof}}. The lemma can be proved in a similar way as Theorem 3.3.1 in \cite{S99}.

\begin{lem}
The transpose $\mathfrak{I}^{\dag}: H^k(S_{-}\p M, \R^{2d}) \to H^k(S(\tau_2M, \R^{2d}))$ of the operator $\mathfrak{I}$ has the following representation:
\[
\mathfrak{I}^{\dag}(h)_{ij}(x)= \int_{S_xM} \Phi^{\dag}(x, \xi) h^{\sharp}(x,\xi) \xi_i\xi_j d\mu_x(\xi).
\]
where $h^{\sharp}$ is the unique lift of $h$ to $SM$ which is invariant under the geodesic flow and
is equal to $h$ on $S_{-}\p M$.
\end{lem}

{\bf{Proof}}. For any $h\in C^{\infty}(S_{-}\p M, \R^{2d})$, we have
\beas
(\mathfrak{I}(\Pi), h) &=& \int_{S_{-}\p M}|\langle \nu(x_0), \xi_0\rangle|d\Sigma^{2d-2}(x_0, \xi_0) h(x_0, \xi_0)\int_{0}^{\infty} \Phi(\mathcal{H}^t(x_0, \xi_0)) \Pi_{ij}(x(t))\xi^i(t)\xi^j(t)dt \\
&=& \int_{SM} \langle h^{\sharp}(x, \xi), \Phi(x, \xi)\Pi_{ij}(x)\xi^i\xi^j \rangle d\mu(x, \xi) \quad \quad (\mbox{by Santalo's formula, see \cite{S99}})\\
&=& \int_{SM} \langle \Phi^{\dag}(x, \xi)h^{\sharp}(x, \xi), \Pi_{ij}(x)\xi^i\xi^j \rangle d\mu(x, \xi) \\
&=& \int_{SM} \sum_{1\leq i,j\leq d} \langle \Phi^{\dag}(x, \xi)h^{\sharp}(x, \xi)\xi^i\xi^j, \Pi_{ij}(x) \rangle d\mu(x, \xi) \\
&=& \left( \int_{S_xM}\Phi^{\dag}(x, \xi)h^{\sharp}(x, \xi)\xi_i\xi_j d\mu_x(\xi), \Pi_{ij}(x) \right),
\eeas
whence the claim follows.

\begin{lem}
 The normal operator $\mathfrak{M}=\mathfrak{I}^{\dag}\mathfrak{I}$ has the following representation:
 \beas
 \mathfrak{M} (\Pi)_{ij}(x)
 &=& \int_{S_xM} d \mu_x(\xi)\int_{-\infty}^{\infty} \Phi^{\dag}(x, \xi) \Phi(\mathcal{H}^t(x,\xi))\,\xi_i \,\xi_j
  \,\Pi_{mn}(\exp_x t\xi)\cdot \dot{\exp}^m_x t\xi \cdot \dot{\exp}^m_x t\xi \, dt \\
 &=& \int_{T_xM} W^{mn}_{ij}(x, \xi) f_{mn}(\exp_{x} \xi) \,d \mu_x(\xi),
  \eeas
 where
\begin{align*}
W^{mn}_{ij}(x, \xi) &= \f{\xi_i \xi_j}{|\xi|^{d+1}} \dot{{\exp}}^m_x (\f{\xi}{|\xi|}) \dot{\exp}^n_x (\f{\xi}{|\xi|}) \\
 &  \cdot\left(  \Phi^{\dag}(x, \f{\xi}{|\xi|}) \Phi(\exp_{x} \xi, \dot{\exp}_x\f{\xi}{|\xi|})
       +  \Phi^{\dag}(x, -\f{\xi}{|\xi|}) \Phi(\exp_x\xi, -\dot{\exp}_x\f{\xi}{|\xi|}) \right).
\end{align*}
\end{lem}



We now show some local properties of the normal operator $\mathfrak{M}$.

Let $\tilde{M}\Subset M$. From now on, we fix $x_*\in \tilde{M}$. We first decompose $\mathfrak{M}$ locally
into two parts based on the separation of singularities of its Schwartz kernel. By the existence of uniformly normal neighborhood
in Riemannian manifold,
we can find $\epsilon_2>0$ and a neighborhood of $x_*$, say $\tilde{U}(x_*)\subset \R^d$,
such that
\be
\exp(x, \cdot)|_{|\xi|< 2 \epsilon_2} \,\, \mbox{is a diffeomorphism for any} \,\,x\in \tilde{U}(x_*).
\ee

Let $\chi_{\epsilon_2} \in C_0^{\infty}(\R)$ be such that $\chi(t)=1$ for $|t|< \epsilon_2$ and
$\chi(t)=0$ for $|t|>2\epsilon_2$. We then define
\bea
(\mathfrak{M_1}f)_{ij}(x)&=& \int_{T_{x}M}W^{mn}_{ij}(x, \xi)f_{mn}(\exp_x\xi))\chi_{\epsilon_2}(|\xi|)\,d\mu_x(\xi),\\
(\mathfrak{M_2}f)_{ij}(x)&=& \int_{T_{x}M}W^{mn}_{ij}(x, \xi)f_{mn}(\exp_x \xi)(1-\chi_{\epsilon_2}(|\xi|))\,d\mu_x(\xi).
\eea
Note that for any $f$ supported in $\tilde{M}$, $f(\exp_x \xi)=0$ for all $|\xi| > T$. Thus we have
$$
(\mathfrak{M_2}f)_{ij}(x)= \int_{\xi\in T_{x}M, \,\, \epsilon_2 < |\xi|< T} W^{mn}_{ij}(x, \xi)f_{mn}(\exp_x\xi)(1-\chi_{\epsilon_2}(|\xi|))\,d\mu_x(\xi).
$$

It is clear that $\mathfrak{M} f= \mathfrak{M_1} f+\mathfrak{M_2} f$.
This gives the promised decomposition of $\mathfrak{M}$.
We next study $\mathfrak{M_1}$ and $\mathfrak{M_2}$
separately.

We first investigate $\mathfrak{M_1}$. In the uniformly normal neighborhood of of $x_*$, by a change of coordinate $y=\exp_x \xi$ and using the fact that $\f{\xi}{|\xi|}= -\nabla_x\rho(x, y)$,
$\dot{\exp}_{x}\xi = \nabla_y\rho(x, y)$, we can deduce that (see \cite{SU04}) $\mathfrak{M_1}$
has the following representation:
\be
(\mathfrak{M_1}f)_{ij}(x)= \f{1}{\sqrt{\det g(x)}}\int_{M}K^{mn}_{ij}(x, y)f_{mn}(y)\chi_{\epsilon_1}(\rho(x,y))\,dy,
\ee
where
\begin{align*}
K^{mn}_{ij}(x, y)= &
\f{\p \rho}{\p x^i}\f{\p \rho}{\p x^j}g^{mm'}(y)g^{nn'}(y) \f{\p \rho}{\p x^{m'}}\f{\p \rho}{\p x^{n'}} \\
& \cdot \left(
\Phi^{\dag}(x, -\nabla_x\rho)\Phi(y, \nabla_y\rho)+\Phi^{\dag}(x, \nabla_x\rho)\Phi(x, -\nabla_y\rho)
\right).
\end{align*}
Following the same argument as in the proof of Lemma 3 in \cite{SU08AJM}, we conclude that the following result holds.
\begin{lem} \label{lem-n11}
 $\mathfrak{M}_1$ is a classic $\Psi$DO of order $-1$ in a neighborhood of $x_*$ with principal symbol
\be
 \sigma_p(\mathfrak{M}_1)^{mn}_{ij}(x, \omega) = \pi \int_{S_xM} (\Phi^{\dag}\Phi)(x, \xi) \xi_i\xi_j\xi^m\xi^n \delta(\omega \cdot \xi) d\mu_x(\xi).
\ee
Here for each fixed $i,j,m,n$, $\sigma_p(\mathfrak{M}_1)^{mn}_{ij}(x, \omega)$ is a matrix from $\R^{2d}$ to $\R^{2d}$. 
\end{lem}


We now proceed to study the operator $\mathfrak{M_2}$ whose property is determined by the exponential map
$\exp(x_*, \cdot)$.
We shall study the operator germ $\mathfrak{M}_{2, \xi_*}$ for each $\xi \in T_{x_*}M$.
We first consider the case when
$\xi_*$ is not a conjugate vector, i.e. $\xi_*$ is a regular vector.

\begin{lem} \label{lem-x-ray5}
Let $\xi_* \in S^*_{x_*}\R^d$ be a regular vector, then there exists a neighborhood $U(x_*)$ of $x_*$ and a neighborhood
$B(x_*, \xi_*)$ of $(x_*, \xi_*)$ such that for any
$\chi\in C^{\infty}_0(B(x_*, \xi_*))$ the following operator
\be
(\mathfrak{M}_{2, \xi_*} f)_{ij}(x)= \int_{T_{x}M} W^{mn}_{ij}(x, \xi)f_{mn}(x,\xi)(1- \chi_{\epsilon_2}(|\xi|))\cdot \chi(x, \xi)\,d\mu_x(\xi) \nonumber
\ee
is a smoothing operator from $\mathcal{E}'(S(\tau_2M), \R^{2d})$ into $C^{\infty}(S(\tau_2\overline{U(x_*)}), \R^{2d})$.
\end{lem}

We next consider the case when $\xi_*$ is a fold vector.
We have the following result.

\begin{lem} \label{lem-modified-n2}
Let $\xi_{*}$ be a fold vector of the map $\exp_{x_*}(\cdot)$.
Then there exists a small neighborhood $U(x_*)$ of $x_*$
and a small neighborhood $B(x_*, \xi_*)$ of $(x_*, \xi_*)$ in $\R^{2d}$
such that for any $\chi \in  C^{\infty}_{0}(B(x_*, \xi_*))$, the operator
$\mathfrak{M}_{2, \xi_*}:  \mathcal{E}'(S(\tau_2M), \R^{2d}) \to  \mathcal{D}'(S(\tau_2U(x_*)), \R^{2d})$ defined by
\be \label{formula-m2}
(\mathfrak{M}_{2, \xi_*} f)_{ij}(x)= \int_{T_{x}M}W^{mn}_{ij}(x, \xi)f_{mn}(x,\xi)(1- \chi_{\epsilon_2}(|\xi|))\cdot \chi(x, \xi) \,d\mu_x(\xi),
\quad f\in \mathcal{E}'(S(\tau_2M), \R^{2d})
\ee
is an FIO of order $-\f{d}{2}$ whose associated canonical relation is compactly supported in the following set
\bea \label{canonical-relation}
\Big\{(x, \xi, y, \eta);  && x\in U(x_*), y=\exp(x, \omega), (x, \omega) \in B(x_*, \xi_*), \, \, \det d_{\omega}\exp(x, \omega)=0,  \nonumber\\
&& \xi= -\eta_i\f{\p \exp^i(x, \omega)}{\p x}, \, \eta \in \mbox{Coker} \,(d_{\omega}\exp(x, \omega)).\Big\}
\eea
Moreover, the canonical relation is the graph of a homogeneous canonical transformation from a neighborhood
$(\exp_{x_*}\xi_*, \dot{\exp}_{x_*}\xi_* ) \in T\R^d$ to $(x_*, \xi_*)\in T\R^d$, and hence
$\mathfrak{M}_{2, \xi_*}$ is bounded from $H_{\m{comp}}^s$ to $H_{\mbox{local}}^{s+d/2}$.
\end{lem}

{\bf{Proof}}. We note that for each fixed index quadruple $(i,j,m,n)$, the Schwartz kernel of the integral operator defined on the right hand side of (\ref{formula-m2}) has the same type of singularities as those discussed in \cite{SUAPDE, BZ12}. Thus a similar argument as therein yields the result.


\section{The X-ray transform operator $\mathfrak{X}$ and its normal $\mathfrak{N}$} \label{sec-x-n}

We study the X-ray transform $\mathfrak{X}$ and its normal $\mathfrak{N}$ in this section. Our goal is to prove Theorem \ref{thm-1} and Corollary \ref{thm-2} and \ref{thm-3}.

We first show that the X-ray transform $\mathfrak{X}$ vanishes on the potential fields.
\begin{lem} \label{lem-ns}
 $ \mathfrak{X}\mathcal{P} = 0$.
\end{lem}
{\bf{Proof}}. It suffices to show that $\mathfrak{X} d^sv=0$ for all $v\in C_0^{\infty}(M)$.
Indeed, let $v\in C_0^{\infty}(M)$. Consider the one-parameter family of diffeomorphisms $\phi^{\tau}: M \to M$ defined by
$$
\f{d \phi^{\tau}}{d \tau} = v, \quad \phi^0= Id.
$$
It is clear that $\phi^{\tau}|_{\p M} =Id$. Denote $g^{\tau}= (\phi^{\tau})^*g$ and $\mathcal{H}(\tau)= \mathcal{H}(g^{\tau})$.
We have $\mathcal{H}(\tau) \equiv \mathcal{H}(0)$. Thus $\mathcal{H}'(0)=0$.
On the other hand, a direct calculation shows that $\f{d g^{\tau} }{d \tau}|_{\tau=0}= d^s v$, which further implies that
$\mathcal{H}'(0)= \mathfrak{X}\f{d g^{\tau} }{d \tau}|_{\tau=0}= \mathfrak{X} d^s v$.
Therefore, we can conclude that
$ \mathfrak{X} d^s v =0$, which completes the proof of the lemma.

\medskip

Note that $\mathfrak{X}= \mathfrak{I}\circ \iota\circ \mathcal{L}$.
We next present a decomposition of its normal $\mathfrak{N}$ which follows from that of the operator $\mathfrak{M}$:
$$
\mathfrak{N}=\mathfrak{X}^{\dag}\mathfrak{X}=
(\iota\circ \mathcal{L})^{\dag}\circ\mathfrak{M}\circ (\iota\circ \mathcal{L}) =
(\iota\circ \mathcal{L})^{\dag}\circ \mathfrak{M}_1 \circ (\iota\circ \mathcal{L}) +
 (\iota\circ \mathcal{L})^{\dag}\circ \mathfrak{M}_2 \circ (\iota\circ \mathcal{L})
 =\mathfrak{N}_1 +\mathfrak{N}_2.
$$

Similar to the truncation of the operator $\mathfrak{X}$, we introduction a truncation for the operator $\mathfrak{I}$.
For any $\alpha\in C^{\infty}_0(S_{-}\p M)$, we define
\begin{equation} \label{formula-I-alpha}
\mathfrak{I}_{\alpha}f(x_0, \xi_0)
= \alpha(x_0, \xi_0) \mathfrak{I}_{\alpha}f(x_0, \xi_0).
\end{equation}
It is clear that
$$\mathfrak{X}_{\alpha}=\mathfrak{I}_{\alpha}\circ \iota\circ \mathcal{L}.$$

By replacing the weight $W^{mn}_{ij}$ with the new one $\alpha^{\sharp}\cdot W^{mn}_{ij}$, where $\alpha^{\sharp}$ is the unique lift of $\alpha$ to $S\R^d$ which is constant along each orbit of the geodesic flow $\mathcal{H}^t(g)$,  we can define
$\mathfrak{M}_{\alpha}$, $\mathfrak{M}_{1, \alpha}$, $\mathfrak{M}_{2, \alpha}$ and consequently $\mathfrak{N}_{\alpha}$, $\mathfrak{N}_{\alpha,1}$ and $\mathfrak{N}_{\alpha,2}$.
It is clear that
\[
\mathfrak{N}_{\alpha}=
(\iota\circ \mathcal{L})^{\dag}\circ\mathfrak{M}_{\alpha}\circ (\iota\circ \mathcal{L}) =
 (\iota\circ \mathcal{L})^{\dag}\circ\mathfrak{M}_{\alpha,1}\circ (\iota\circ \mathcal{L})+  (\iota\circ \mathcal{L})^{\dag}\circ\mathfrak{M}_{\alpha,2}\circ (\iota\circ \mathcal{L})= \mathfrak{N}_{\alpha,1}+\mathfrak{N}_{\alpha,2}.
\]

\subsection {Local properties of the normal operator $\mathfrak{N}$}

Let $x\in \tilde{M}\Subset M$ be a fold-regular point with the compact subset $\mathcal{Z}(x_*)\subset S_{x_*}M$ in Definition \ref{def-fold-regular-point}. We now construct
a cut-off function $\alpha \in C_0^{\infty}(S_{-}\p M)$ such that
in a neighborhood of $x_*$, $\mathfrak{N}_{\alpha, 1}$ is elliptic and $\mathfrak{N}_{\alpha, 2}$ is smoother than
$\mathfrak{N}_{\alpha, 1}$ (see Lemma \ref{lem-n1} and \ref{lem-n2}).
The idea is to select a complete set of geodesics with no conjugate points expect strong fold-regular ones. We follow closely the argument in \cite{BZ12}.

Denote
$\mathcal{C}_{\epsilon_2, T}\mathcal{Z}= \{r\xi; \xi\in \mathcal{Z}(x_*), r\in \R \quad \mbox{and} \,\,\, \epsilon_2\leq |r|\leq T\}$. With the help of
Lemma \ref{lem-x-ray5} and Lemma \ref{lem-modified-n2}, there exist a finite number of vectors
$\xi_j \in \mathcal{C}_{\epsilon_2, T}\mathcal{Z}$, $j=1, 2, ... N$ such that for each $\xi_j$,
there exist two neighborhoods $B_0(x_*, \xi_*) \Subset B(x_*, \xi_*)$ of $(x_*, \xi_j)\in \R^{2d}$
and a function $\chi_j \in C^{\infty}_0(B(x_*, \xi_j))$ with the following conditions:

(1). The operator $\mathfrak{M}_{2, \xi_j}$ defined by
\[
(\mathfrak{M}_{2, \xi_j} f)_{ij}(x)= \int_{T_{x}M}W^{mn}_{ij}(x, \xi)f_{mn}(\exp_x \xi)(1- \chi_{\epsilon_2}(|\xi|))\cdot \chi_j(x,\xi)\,d\mu_x(\xi)
\]
is bounded from $H^k(S(\tau_2(M, M), \R^{2d}))$ to $H^{k+\f{d}{2}}(S(\tau_2U(x_*, \xi_*)), \R^{2d})$.

(2). $\mathcal{C}_{\epsilon_2, T} \mathcal{Z} \subset \bigcup_{j=1}^{N} B_0(x_*, \xi_j)$.

Denote by $\mathcal{A}_0$ be the greatest connected open symmetric subset in $\bigcup_{j=1}^{N} B_0(x_*, \xi_j)$
which contains $\mathcal{C}_{\epsilon_2, T} \mathcal{Z} $. Here and after, we say that a set $\mathcal{B}$
in $\R^{2d}$ is symmetric if $(x, \xi)\in \mathcal{B}$ implies that $(x, -\xi)\in \mathcal{B}$.
Define
$$
\mathcal{A}_{\epsilon} =\{(x, \xi) \in \R^{2d}: |x-x_*|\leq \epsilon, \,\, \epsilon_2 \leq |\xi| \leq T\}
$$
for each $\epsilon>0$. It is clear that $\mathcal{A}_{\epsilon}$ is compact in $\R^{2d}$,
so is the set $\mathcal{A}_{\epsilon}\backslash \mathcal{A}_0$.

\begin{lem} \label{lem-cutoff}
There exist $\epsilon_3>0$ and $\alpha \in C_0^{\infty}(S_{-}\p M)$ such that the following two conditions are satisfied:
\bea
\alpha(x_0, \xi_0)&=& 1\quad \mbox{for all} \,\, (x_0, \xi_0) \in \mathcal{H}^t(g)(\mathcal{Z}(x_*)), \label{alpha1}\\
\alpha(x_0, \xi_0)&=& 0 \quad \mbox{for all } \, (x_0, \xi_0)\in \mathcal{H}^t(g)(\mathcal{A}_{\epsilon_3}\backslash \mathcal{A}_0)\label{alpha2}.
\eea

\end{lem}

{\bf{Proof}}. See \cite{BZ12}.

\begin{lem} \label{lem-m1}
There exists a neighborhood $U(x_*)$ of $x_*$ such that
\be \label{kernel}
 \ker \sigma_p(\mathfrak{M}_{\alpha, 1})(x, \omega) = \{v_i^k\omega_j+ v_j^k\omega_i:  \,\, \, v_i^k\in \R, 1\leq i, j \leq d, 1\leq k \leq 2d \},
 \ee
 for all $x\in U(x_*)$ and $\omega \in T^*_{x_*}M \backslash 0$.
\end{lem}

{\bf{Proof}}. Let $f\in \ker \sigma_p(\mathfrak{M}_{\alpha, 1})(x, \omega)$, then
\be \label{identity1}
\langle   f, \sigma_p(\mathfrak{M}_{\alpha, 1})(x, \omega)f \rangle =0.
\ee
Similar to Lemma \ref{lem-n11}, the symbol of $\mathfrak{M}_{\alpha,1}$ has the following form:
\be \label{identity2}
 \sigma_p(\mathfrak{M}_{\alpha,1})^{mn}_{ij}(x, \omega) = \pi \int_{\xi \in S_xM} |\alpha^{\sharp}(x, \xi)|^2(\Phi^{\dag}\circ\Phi)(x, \xi) \xi_i\xi_j\xi^m\xi^n \delta(\omega \cdot \xi) d\mu_x(\xi).
\ee
Substituting (\ref{identity2}) into (\ref{identity1}), we can deduce that
$$
\int_{\xi\in S_xM, \,\,  \xi \perp \omega} |\alpha^{\sharp}(x, \xi)|^2 \langle  \Phi(x, \xi)f_{ij}\xi^i\xi^j, \Phi(x, \xi)f_{ij}\xi^i\xi^j \rangle d\mu_x(\xi)=0.
$$
We now show that (\ref{kernel}) holds for the point $x_*$.
Indeed, by Condition (\ref{alpha1}) and the ``completeness'' property of the set $\mathcal{Z}(x_*)$, there exists $\xi_0 \in S_xM$ with $\xi_0 \perp \omega$ such that $\alpha^{\sharp}(x, \xi_0)=1$.
It follows that
$$
\langle  \Phi(x, \xi)f_{ij}\xi^i\xi^j, \Phi(x, \xi)f_{ij}\xi^i\xi^j \rangle=0
$$
for $\xi$ in a small neighborhood of $\xi_0$ such that $\xi \in S_xM$ and $\xi \perp \omega$.
Since $\Phi(x, \xi)$ is invertible, this further implies that
$$
f_{ij}\xi^i\xi^j=0, \quad \mbox{for $\xi$ in a small neighborhood of $\xi_0$ such that $\xi \in S_xM$ and $\xi \perp \omega$}.
$$

Note that for each $k=1, 2,...2d$, the following decomposition holds for $f^k=\{f^k_{ij}; 1\leq i, j\leq d\}$:
$$
 f^k = h^k + i_{\omega} v^k, \quad \mbox{with} \, \, j_{\omega} h^k =0.
$$
Thus,
$\left(h^k_{ij} + (i_{\omega} v^k)_{ij}\right)\xi^i\xi^j=0$ for $\xi$ in a small neighborhood of $\xi_0$ such that $\xi \in S_xM$ and $\xi \perp \omega$. This together with the fact that $(i_{\omega} v^k))_{ij}\xi^i\xi^j=0$ yield $h^k_{ij}\xi^i\xi^j=0$ for $\xi \perp \omega$. On the other hand, the equality $j_{\omega} h^k=0$ implies that
$h^k_{ij} \omega^i=0$. Therefore, we can conclude that $h^k_{ij}\xi^i\xi^j=0$ for all $\xi$ in a small neighborhood of $\xi_0$. It follows that
$h^k =0$. As a result, we get $f^k= i_{\omega} v^k = \f{1}{2}(v_i^k\omega_j+ v_j^k\omega_i )$.
This proves the lemma for the point $x_*$.

To prove (\ref{kernel}) for other points, we exploit the continuity of the function $\alpha^{\sharp}$. In fact, we can find a neighborhood $U(x_*)$ of $x_*$ such that for each $x\in U(x_*)$, there exists a ``complete'' set $\mathcal{Z}(x) \in S_xM$ such that $\alpha^{\sharp}(x, \xi)> \f{1}{2}$ for all $\xi \in \mathcal{Z}(x)$. Then a similar argument as for the point $x_*$ proves that (\ref{kernel}) holds for these points in $U(x_*)$. This completes the proof of the lemma.

\begin{lem} \label{lem-n111}
$\mathfrak{N}_{\alpha, 1}$ is a $\Psi$DO of order one and there exists a neighborhood $U(x_*)$ of $x_*$ such that for $x \in U(x_*)$,
\[
 \ker \sigma_p(\mathfrak{N}_{\alpha, 1})(x, \omega)= \{f\in S_2: f_{ij}=v_i\omega_j+ v_j\omega_i \,\, \mbox{for some}\,\, v\in S_1\}.
\]
Moreover,  $\sigma_p(\mathfrak{N}_{\alpha, 1})(x, \omega)$ is positive definite on the set $\{ f \in S_2:  f_{ij} \omega^i =0\}$.
\end{lem}

{\bf{Proof}}. We first recall that
$\mathfrak{N}_{\alpha, 1}= (\iota \circ \mathcal{L})^{\dag} \circ \mathfrak{M}_{\alpha, 1} \circ (\iota \circ \mathcal{L})$. Since
$\mathfrak{M}_{\alpha, 1}$ is a $\Psi$DO of order $-1$ and both $\iota \circ \mathcal{L}$ and $(\iota \circ \mathcal{L})^{\dag}$
are differential operators of order one, $\mathfrak{N}_{\alpha, 1}$ is a $\Psi$DO of order one.

We next determine the kernel of $\sigma_p(\mathfrak{N}_{\alpha, 1})$. Let $f=\{f_{ij}\} \in \ker \sigma_p(\mathfrak{N}_{\alpha, 1})(x, \omega)$.
Note that
$$ \sigma_p(\mathfrak{N}_{\alpha, 1})= \sigma_p((\iota \circ \mathcal{L})^{\dag})\sigma_p(\mathfrak{M}_{\alpha, 1})\sigma_p(\iota \circ \mathcal{L})=
(\sigma_p(\iota \circ \mathcal{L}))^{\dag}\sigma_p(\mathfrak{M}_{\alpha, 1})\sigma_p(\iota \circ \mathcal{L}).$$
We have
$$
\sigma_p(\iota \circ \mathcal{L})(x, \omega)f  \in \ker \sigma_p(\mathfrak{M}_{\alpha, 1})(x, \omega).
$$
By Lemma \ref{lem-m1}, there exists $v=\{v_{i}^k; 1\leq k \leq 2d, 1 \leq i \leq d\}$ such that
$$
(\sigma_p(\iota \circ \mathcal{L})(x, \omega)f )_{mn}^k = \f{1}{2}(v_m^k\omega_n + v_n^k \omega_m), \quad k=1, 2,...2d.
$$
On the other hand, a direct calculation shows that
$$
(\sigma_p(\iota \circ \mathcal{L})(x, \omega)f )_{mn}^{k+d} = \f{1}{2}g^{kl}( f_{ml}\omega_n + f_{nl}\omega_m - f_{mn}\omega_l ), \quad k=1, 2,...d.
$$
Therefore,
\be \label{equations1}
 f_{ml}\omega_n + f_{nl}\omega_m - f_{mn}\omega_l = v_{m l}\omega_n + v_{nl}\omega_m, \quad l=1, 2,...d.
\ee
where $v_{mn} := g_{ln}(x)v^{l+d}_{m}$, $1\leq m, n, l \leq d$.

We now solve (\ref{equations1}) for $f$. Let $\omega_a$ be the component of $\omega$ such that $\omega_a \neq 0$.
First, by setting $m=n=a$ in (\ref{equations1}), we obtain
\be \label{equations2}
 2f_{al}\omega_a - f_{aa}\omega_l = 2v_{al}\omega_a.
\ee
Next, by setting $l=a$ in (\ref{equations2}), we further get
$$
f_{aa} = 2v_{aa}.
$$
Substituting this back into (\ref{equations2}), it follows that
\be \label{equations3}
f_{al}= \f{\omega_l}{\omega_a} v_{aa} + v_{al}.
\ee
Finally, substituting (\ref{equations3}) into (\ref{equations1}) and letting $l=a$, we deduce that
\beas
f_{mn} &=& \f{2v_{aa}}{\omega_a^2} \omega_m\omega_n + \f{\omega_m}{\omega_a}(v_{an}-v_{na}) +\f{\omega_n}{\omega_a}(v_{am}-v_{ma}) \\
 &=&  \f{\omega_m}{\omega_a}(v_{an}-v_{na}+ \f{v_{aa}}{\omega_a} \omega_n) + \f{\omega_n}{\omega_a}(v_{am}-v_{ma}+ \f{v_{aa}}{\omega_a} \omega_m) \\
 &=&  \omega_m u_n + \omega_n u_m,
\eeas
where $u_m = \f{1}{\omega_a}(v_{am}-v_{ma})+ \f{v_{aa}}{\omega_a^2} \omega_m$.
This gives the desired solution $f$.
From this, we conclude that
\[
 \ker \sigma_p(\mathfrak{N}_{\alpha, 1})(x, \omega) = \{f \in S_2: f_{ij}=v_i\omega_j+ v_j\omega_i \},
\]
for $x\in U(x_*)$ where $U(x_*)$ is chosen as in Lemma \ref{lem-m1}.
Finally, the remaining part of the lemma follows from the fact that $\sigma_p(\mathfrak{N}_{\alpha, 1})(x, \omega)$ is symmetric and is non-negative.
This completes the proof of the Lemma.
\medskip

\begin{lem} \label{lem-n1}
Let $x_*\in \tilde{M}\Subset M$ be a strong fold-regular point, $\alpha\in C^{\infty}_0(S_{-}\p M)$ be chosen as in Lemma \ref{lem-cutoff} and let $k$ is a nonnegative integer. Then
there exist two neighborhoods $U(x_*)\Subset \tilde{U}(x_*)$ of $x_*$ such that the following estimate holds for all $f \in H^k(S(\tau_2(\tilde{M}, M)))$
\be \label{equ-110}
\|f^s\|_{H^k(U(x_*))} \lesssim  \|\mathfrak{N}_{1, \alpha} f\|_{H^{k-1}(\tilde{U}(x_*))} + \|f\|_{H^{k-1}(M)}.
\ee

\end{lem}

{\bf{Proof}}. We divide the proof into the following three steps.

Step 1. We construct pseudo-inverse of the operator $\mathfrak{N}_{\alpha, 1}$ for solenoidal 2-tensor fields. Consider the operator $\mathcal{A}= |D|^{-1}\mathfrak{N}_{\alpha, 1} + \mathcal{P}$, where $|D|^{-1}$ is a properly supported parametrix for $(-\Delta_g)^{-\f{1}{2}}$ in $M$.  By Lemma \ref{lem-n111}, $\mathcal{A}$
is an elliptic $\Psi$DO of order zero in a neighborhood $\tilde{U}(x_*)$ of $x_*$. We thus can find a $\Psi$DO of order zero, denoted by $\mathcal{B}$, such that
$$
\mathcal{B}\mathcal{A} = \mathcal{Q} + \mathcal{K},
$$
where $\mathcal{Q}$ is a $\Psi$DO with symbol $\sigma(\mathcal{Q})(x, \omega)= Id$ for $x\in \tilde{U}(x_*)$ and $\mathcal{K}$
is a smoothing operator. Let $\mathcal{K}_1 = \mathcal{K} + \mathcal{Q}- Id$. Then $\mathcal{B}\mathcal{A}= Id +\mathcal{K}_1 $, i,e.
$$
\mathcal{B}|D|^{-1}\mathfrak{N}_{\alpha, 1} + \mathcal{B}\mathcal{P} = Id +\mathcal{K}_1.
$$
Applying $\mathcal{S}$ from the right on both sides and using the fact that $\mathcal{P}\mathcal{S}=0$, we obtain
$$
\mathcal{B}|D|^{-1}\mathfrak{N}_{\alpha, 1}\mathcal{S} = \mathcal{S} +\mathcal{K}_1\mathcal{S}.
$$
This completes the construction of the pseudo-inverse of $\mathfrak{N}_{\alpha, 1}$.
\medskip

Step 2. Let $f\in H^k(S(\tau_2(\tilde{M}, M)))$,
we show that
\be  \label{equ-111}
\|f^s\|_{H^k(U(x_*))}
 \lesssim  \|\mathfrak{N}_{1, \alpha} f^s\|_{H^{k-1}(\tilde{U}(x_*))} + \|f\|_{H^{-t}(M)}
\ee
for any $t>0$. We argue as follows. By the result in Step 1, we have
$$
f^s = \mathcal{S}f^s =\mathcal{B}|D|^{-1}\mathfrak{N}_{\alpha, 1}f^s - \mathcal{K}_1f^s \quad \mbox{in} \,\, \tilde{U}(x_*).
$$
Let $U(x_*)$ be another neighborhood of $x_*$ such that $U(x_*) \Subset \tilde{U}(x_*)$. Then
\[
\|f^s\|_{H^k(U(x_*))}  \leq   \|\mathcal{B}|D|^{-1}\mathfrak{N}_{\alpha, 1}f^s\|_{H^k(U(x_*))} +  \|\mathcal{K}_1f^s\|_{H^k(U(x_*))}.
\]
We now estimate the term $\|\mathcal{B}|D|^{-1}\mathfrak{N}_{\alpha, 1}f^s\|_{H^k(U(x_*))}$.
Let $\chi_{U(x_*)}$ be a smooth function such that $\chi_{U(x_*)}|_{U_1(x_*)}=1$ for some $U_1(x_*)$ satisfying
$U(x_*) \Subset U_1(x_*) \Subset \tilde{U}(x_*)$, and $\chi_{U(x_*)}=0$ outside $\tilde{U}(x_*)$.
Note that
$$
\mathcal{B}|D|^{-1}\mathfrak{N}_{\alpha, 1}f^s= \mathcal{B}|D|^{-1} (\chi_{\tilde{U}(x_*)}\mathfrak{N}_{1, \alpha} f^s)+
\mathcal{B}|D|^{-1} (1-\chi_{\tilde{U}(x_*)})\mathfrak{N}_{1, \alpha} f^s.
$$
We can deduce that
\begin{align*}
&\|\mathcal{B}|D|^{-1} (\chi_{U(x_*)}\mathfrak{N}_{1, \alpha}f^s)\|_{H^k(U(x_*))}  \lesssim  \|\chi_{U(x_*)}\mathfrak{N}_{1, \alpha} f^s\|_{H^{k-1}(\tilde{U}(x_*))}\lesssim \|\mathfrak{N}_{1, \alpha} f^s\|_{H^{k-1}(\tilde{U}(x_*))}, \\
& \|\mathcal{B}|D|^{-1} (1-\chi_{U(x_*)})\mathfrak{N}_{1, \alpha} f^s\|_{H^k(U(x_*))}  \lesssim \|f\|_{H^{-t}(M)},
\end{align*}
where for the second inequality above we used the fact that the operator $\mathcal{B}|D|^{-1} (1-\chi_{U(x_*)})\mathfrak{N}_{1, \alpha}\mathcal{S}$ is smoothing from $H^k(S(\tau_2(\tilde{M}, M)))$ to $C^{\infty}(\overline{S(\tau_2U(x_*))})$. Thus, we conclude that
\[
\|\mathcal{B}|D|^{-1}\mathfrak{N}_{\alpha, 1}f^s\|_{H^k(U(x_*))} \lesssim \|\mathfrak{N}_{\alpha, 1}f^s\|_{H^{k-1}(U(x_*))} + \|f\|_{H^{-t}(M)}.
\]
for any $t>0$.

To finish the proof of (\ref{equ-111}), it remains to show that
$\mathcal{K}_1\mathcal{S}$ is a smoothing operator from $H^k(S(\tau_2(\tilde{M}, M)))$ to $C^{\infty}(S(\tau_2\tilde{U}(x_*)))$.
Indeed, note that $\mathcal{K}_1 = \mathcal{K} + \mathcal{Q}- Id$. We need only show that $(\mathcal{Q}- Id)\mathcal{S}$ is smoothing from
$L^2(S(\tau_2(\tilde{M}, M)))$ to $C^{\infty}(S(\tau_2\tilde{U}(x_*)))$. But this is a consequence of the fact that the symbol of
$\mathcal{Q}- Id$ vanishes for $x\in \tilde{U}(x_*)$.

\medskip

Step 3. We finally prove (\ref{equ-110}). Observe that $\mathfrak{N}_{\alpha, 1}f^s= \mathfrak{N}_{\alpha, 1}f - \mathfrak{N}_{\alpha, 1}\mathcal{P}f$.
By Lemma \ref{lem-n111}, the principle symbol of the operator $\mathfrak{N}_{\alpha, 1}\mathcal{P}$ is equal to zero for $x\in \tilde{U}(x_*)$. As a result, $\mathfrak{N}_{\alpha, 1}\mathcal{P}$ is a $\Psi$DO of order zero and therefore the estimate below holds
\[
 \|\mathfrak{N}_{\alpha, 1}\mathcal{P}f\|_{H^{k-1}(U(x_*))} \lesssim \|f\|_{H^{k-1}(M)}.
\]
It follows that
\[
\|\mathfrak{N}_{\alpha, 1}f^s\|_{H^{k-1}(U(x_*))} \lesssim \|\mathfrak{N}_{\alpha, 1}f\|_{H^{k-1}(U(x_*))}+ \|f\|_{H^{k-1}(M)}.
\]
By substituting the above inequality into (\ref{equ-111}), we get (\ref{equ-110}).
This completes the proof of the Lemma.

\begin{lem} \label{lem-n2}
Let $x_*\in \tilde{M}\Subset M$ be a strong fold-regular point and $\alpha\in C^{\infty}_0(S_{-}\p M)$ be chosen as in Lemma \ref{lem-cutoff}. Then
there exists a neighborhood $U(x_*)$ of $x_*$ such that the following estimate holds for all
$f\in H^k(S^2(\tau_2(\tilde{M}, M)))$
\[
\|\mathfrak{N}_{2, \alpha} f\|_{H^{k-1}(U(x_*))}  \lesssim  \|f\|_{H^{k-\f{d}{2}+1}(M)}
\]
where $k$ is any nonnegative integer.
\end{lem}

{\bf{Proof}}: With the help of Lemma \ref{lem-modified-n2}, the proof is similar to that of Lemma 6.7 in \cite{BZ12}.

\subsection{Proof of Theorem \ref{thm-1} and Corollary \ref{thm-2} and \ref{thm-3}} \label{sub-sec-proof-x-ray}

{\bf{Proof of Theorem \ref{thm-1}}}: It is a direct consequence of Lemma \ref{lem-n1} and \ref{lem-n2}.

\medskip

\noindent{\bf{Proof of Corollary \ref{thm-2}}}: For each $x\in \tilde{M}$, by Theorem \ref{thm-1}, there exist neighborhoods $U(x) \Subset \tilde{U}(x_*)$ of
$x$ and a smooth function $\alpha \in C^{\infty}_0(S_{-}\p M)$ such that the estimate (\ref{estimate1}) holds. Since $\tilde{M}$ is compact, we can find finite number of points, say $x_1$, $x_2$, ... $x_N$, such that
$\tilde{M} \subset \bigcup_{j=1}^{N}U(x_j)$ and the following estimate holds for each $j$:
\[
 \|f^s\|_{H^k(U(x_j))} \lesssim \|\mathfrak{N}_{\alpha_j}f\|_{H^{k-1}(\tilde{U}(x_j))} +  \|f\|_{H^{k+1-\f{d}{2}}(M)}.
\]

Let $\tilde{M}_1= \bigcup_{j=1}^{N}U(x_j)$. Then $\tilde{M} \Subset \tilde{M}_1$. We claim that
$$
 \|f^s\|_{H^k(M \backslash \tilde{M}_1)}  \lesssim   \|f\|_{H^{k+1-\f{d}{2}}(M)}.
$$
Indeed, note that $f^s = \mathcal{S}f = (Id- d^s(\triangle_M^s)^{-1}\delta^s) f$, where $\triangle_M^s$ is the Dirichlet realization of the operator
$\triangle^s := \delta^s d^s$ in $M$ and $(\triangle_M^s)^{-1}$ is its inverse. Using the pseudo-local property of the operator $\mathcal{S}$, we see that $\mathcal{S}$ is smoothing from $L^2(\tilde{M})$ to $C^{\infty}(\overline{M \backslash \tilde{M}_1})$. Hence the claim follows.

Therefore, we deduce that
\beas
\|f^s\|_{H^k(M)} & \leq & \|f^s\|_{H^k(M\backslash \tilde{M}_1)}  + \sum_{j=1}^N \|f\|_{H^k(U(x_j))}   \\
& \lesssim & \sum_{j=1}^N\|\mathfrak{N}_{\alpha_j}f\|_{H^k(\tilde{U}(x_j))} +  \|f\|_{H^{k+1-\f{d}{2}}(M)}.
\eeas
The corollary is proved.

\medskip
\noindent{\bf{Proof of Corollary \ref{thm-3}}}:

Step 1. Let $x_j$, $U(x_j)$ and $\alpha_j$ be chosen as in Theorem \ref{thm-2}.
Denote by $H$ the Hilbert space $\prod_{j=1}^{N}H^{-1}(S(\tau_2\tilde{U}(x_j)))$.
We consider the operator
$$
T: L^2(S(\tau_2(\tilde{M}, M))) \to H
$$ defined by
$$
T f = (\mathfrak{N}_{\alpha_1}f, \mathfrak{N}_{\alpha_2}f,..., \mathfrak{N}_{\alpha_N}f).
$$
Then Theorem \ref{thm-2} implies that (using the fact that $d\geq 3$)
\be \label{equation61}
\|f^s\|_{L^2(M)}  \lesssim \|Tf\|_H + \|f\|_{H^{-\f{1}{2}}(M)}.
\ee

Note that the operator $T$ vanishes on the potential vector fields (see Lemma \ref{lem-ns}). We can define
$T_0:\mathcal{S}L^2(S(\tau_2(\tilde{M}, M))) \to H$ by
\[
T_0f^s=Tf.
\]

Step 2.
We show that the space
$\mathfrak{L}:=\ker{T_0}$ is finite dimensional in $\mathcal{S}L^2(S(\tau_2(\tilde{M}, M)))$.
Indeed, by contradiction, assume that there exist a pairwise orthogonal sequence $f^s_n \in \mathfrak{L}$ such that
$\|f^s_n\|_{L^2(M)}=1$ and $T_0f^s_n=0$.
By Lemma \ref{lem-aux1}, there exist $C>0$, $g_n \in L^2(S(\tau_2(\tilde{M}, M)))$ such that
\[
\|g_n\|_{L^2(M)} \leq C, \quad \mbox{and }\, g_n^s=f_n^s.
\]
Then the inequality (\ref{equation61}) yields that
\[
\|f^s_n-f^s_m\|_{L^2(M)}  \lesssim  \|g_n-g_m\|_{H^{1-\f{d}{2}}(M)}.
\]

Since $H^{-\f{1}{2}}(M)$ is compactly embedded into $L^2(M)$, we can find a subsequence of $g_n$, still denoted by $g_n$, such that $g_n$ is Cauchy in $H^{-\f{1}{2}}(M)$. Then $f^s_n$ is also Cauchy in $L^2(M)$. This contradicts to the assumption that $f^s_n$ are pairwise orthogonal in $L^2(S(\tau_2(\tilde{M}, M)))$. This contradiction proves our assertion.

\medskip

Step 3.
We finally show (\ref{equation41}).
Assume the contrary, then there exists a sequence $f_n \in L^2(S(\tau_2(\tilde{M}, M)))$ such that
$\{f_n^s\}_{n=1}^{\infty} \subset  \mathcal{S}L^2(S(\tau_2(\tilde{M}, M))) \bigcap \mathfrak{L}^{\perp} $,
$\|f_n^s\|_{L^2(M)} =1$ and  $\|Tf_n\|_{H}\leq \f{1}{n}$ for all $n$.
By the same argument as in Step 2, we may assume that the sequence $\{f_n\}_{n=1}^{\infty}$ is
Cauchy in $L^2(S\tau_2M)$.
As a result, we can conclude that $\{f_n^s\}_{n=1}^{\infty}$ is Cauchy in
$L^2(M)$ by using (\ref{equation61}).
Note that $\mathcal{S}L^2(S(\tau_2(\tilde{M}, M)))$ is closed in $L^2(S(\tau_2M)$ by Lemma \ref{lem-aux2}. There exists $f_0 \in L^2(S(\tau_2(\tilde{M}, M)))$ such that
$f_0^s= \lim_{n\to \infty}f_n^s$. It is clear that
$f_0^s \in \mathfrak{L}$.
However, since $\mathfrak{L}^{\perp}$ is closed, as the limit of a sequence of functions in $\mathfrak{L}^{\perp}$,
$f_0^s$ should also belong to
$\mathfrak{L}_0^{\perp}$. Therefore,  $f_0^s=0$.
But this contradicts to the fact that
$\|f_0^s\|_{L^2(M)}= \lim_{{n\to \infty}} \|f_n^s\|_{L^2(M)}=1$.
This contradiction completes the proof of (\ref{equation41}) and hence the corollary.

\medskip
\noindent{\bf{Proof of Corollary \ref{thm-4}}}: The same argument as for Corollary \ref{thm-3}.




\medskip

We now present two axillary lemmas which are needed in the proofs above.

\begin{lem} \label{lem-aux1}
Let $k$ be a nonnegative integer. Then there exists $C>0$ such that for any $f^s \in \mathcal{S}H^k(S(\tau_2(\tilde{M}, M)))$, we can find
$g\in H^k(S(\tau_2(\tilde{M}, M)))$ such that
\[
\|g\|_{H^k(M)} \leq C \|f^s\|_{H^k(M)}, \quad \mbox{and }\, g^s=f^s.
\]
\end{lem}

{\bf{Proof}}. We first prove the lemma for the case $k=0$.
Let $f=f^s+ d^sv$ where $f\in L^2(S\tau_2(\tilde{M}, M))$ and $v\in H^1_0(S(\tau_2(\tilde{M}, M)))$.
We need find $g\in L^2(S(\tau_2(\tilde{M}, M)))$ such that $g^s=f^s$ and $\|g\|_{L^2(M)} \leq C \|f^s\|_{L^2(M)}$.
We argue as follows.

We first claim that
\be \label{inequality-1}
\|v\|_{H^1(M\backslash \tilde{M})} \lesssim \|f^s\|_{L^2(M)}.
\ee
Indeed, noting that $f=0$ in $M\backslash \tilde{M}$, we have
\[
d^sv=-f^s  \quad \mbox{in}\,\, M\backslash \tilde{M}.
\]
Using $v|_{\p M} =0$ and the same argument as in \cite{SU04}, we can show that
\be \label{inequality-a}
\|v\|_{L^2(M\backslash \tilde{M})} \lesssim \|f^s\|_{L^2(M\backslash \tilde{M})}.
\ee
On the other hand, observe that
\beas
(d^sv)_{ij}&=& \f{1}{2}(\nabla_i v_j + \nabla_j v_i)\\
&=& \f{1}{2}( \f{\p v_j}{\p x^i} + \f{\p v_i}{\p x^j} - 2\Gamma^k_{ij}v_k).
\eeas
Thus,
\[
\f{\p v_j}{\p x^i} + \f{\p v_i}{\p x^j} = 2 (d^sv)_{ij} + 2 \Gamma^k_{ij}v_k = -2f^s +2 \Gamma^k_{ij}v_k \quad \mbox{in}\,\,M\backslash \tilde{M}.
\]
Consequently,
\be \label{inequality-b}
\|\f{\p v_j}{\p x^i} + \f{\p v_i}{\p x^j}\|_{L^2(M\backslash \tilde{M})} \lesssim \|f^s\|_{L^2(M\backslash \tilde{M}))}+
\|v\|_{L^2(M\backslash \tilde{M}))}  \lesssim \|f^s\|_{L^2(M\backslash \tilde{M}))}.
\ee
By using the following Korn's inequality (see for instance \cite{DL76})
\be  \label{inequality-c}
\|v\|_{H^1(M\backslash \tilde{M})} \lesssim
\|v\|_{L^2(M\backslash \tilde{M})} +
 \sum_{ij} \|\f{\p v_j}{\p x^i} + \f{\p v_i}{\p x^j}\|_{L^2(M\backslash \tilde{M})},
\ee
the estimate (\ref{inequality-1}) follows immediately from (\ref{inequality-a})-(\ref{inequality-c}).
This completes the proof of the claim.

Next, for each $v\in H^1(M\backslash \tilde{M})$, we can find an extension $\tilde{v} \in H^1(M)$ such that
$\tilde{v}=v$ in $M\backslash \tilde{M}$ and $\|\tilde{v}\|_{H^1(M)} \lesssim \|v\|_{H^1(M\backslash \tilde{M})}$.
Let $g = f^s + \tilde{v}$. Then $g \in L^2(S(\tau_2(\tilde{M}, M)))$ and $g^s=f^s$.  Moreover, $\|g\|_{L^2(M)} \lesssim \|f^s\|_{L^2(M)}$. The lemma is proved for the case $k=0$.

Finally, for the case $k\geq 1$, the argument is almost the same. The key point is to derive the following estimate
$$
 \|v\|_{H^{k+1}(M\backslash \tilde{M})} \lesssim \|d^sv\|_{H^k(M\backslash \tilde{M}))},
$$
by Korn's inequality and the estimate (\ref{inequality-a}).

\medskip

As a direct consequence of the above lemma, we obtain the following result.

\begin{lem} \label{lem-aux2}
The linear space $\mathcal{S}H^k(S(\tau_2(\tilde{M}, M)))$ is closed in $H^k(S(\tau_2M))$.
\end{lem}

\section{Proof of Theorem \ref{thm-lens}} \label{sec-final}
We write $\tilde{g}= g +f$. Then $\mbox{supp}f \subset \hat{M}$.
We assume that $f$ is sufficiently small in $C^{k,s}$ norm for some sufficiently large integer $k$ and $0< s<1$ (we will choose $k=15$ and $s=\f{1}{2}$ at the end of the proof). We divide the argument into the following eight steps.

\medskip

Step 1.
We first construct a diffeomorphism $\phi$ of $M$ such that $\phi^*\tilde{g}$ is solenoidal with respect to the metric $g$, following the argument in \cite{CDS00}. In fact, by solving the equation $\delta^s(\phi^*\tilde{g}-g)=0$, we obtain a unique $\phi$ which is close to the identity map and equal to it on the boundary $\p M$ and satisfies the following estimate
\be
\|\phi-Id\|_{C^{k, s}} \lesssim \|f\|_{C^{k, s}}.
\ee

We denote $\tilde{g}_1 = \phi^*\tilde{g}$ and $f_1= \tilde{g}_1-g$. Then the following estimate holds
\be \label{inequality0}
\|f_1\|_{C^{k-1, s}} \lesssim \|f\|_{C^{k, s}}.
\ee

Step 2. We construct another diffeomorphism $\psi$ of $M$ such that $\psi^*\tilde{g}_1=e$ in a small neighborhood of $\p M$.  We follow the same approach as in Section 4.2 in \cite{SU09}.
Denote by $\exp_{\p M, g}$ and $\exp_{\p M, \tilde{g}_1}$ the boundary normal coordinate near $\p M$ with respect to the metric $g$ and $\tilde{g}_1$, respectively. Both maps are well-defined from $\p M \times [0, \epsilon_0]$ for some $\epsilon_0$, which can be chosen to be independent of $\tilde{g}$, to some neighborhoods of $\p M$.
Let $\psi_1= \exp_{\p M, g}(\exp_{\p M, \tilde{g}_1})^{-1}$. Then $\psi_1$ maps a small neighborhood of $\p M$ to another and
satisfies the following estimate:
\be
\|\psi_1 -Id\|_{C^{k-2}} \lesssim \|f_1\|_{C^{k-1}}.
\ee
If $\|f_1\|_{C^{k-1}}$ is sufficiently small, we can show that there exists a unique vector field $v\in C^{k-2}$ in a sufficiently small neighborhood of $\p M$, say $W_1$,  with $v=0$ on $\p M$ such that
$$
\psi_1(x)= \exp_xv(x)
$$
and
\[
\|v\|_{C^{k-2}} \lesssim \|\psi_1 -Id\|_{C^{k-2}} \lesssim \|f_1\|_{C^{k-1}}.
\]

We now extend $\psi_1$ to a diffeomorphism of the whole domain $M$.
Let $W$ be another neighborhood of $\p M$ such that $W \Subset W_1\bigcap (M\backslash \hat{M})$. We may choose $W$ such that $M \backslash \overline{W}$ is a strictly convex and smooth sub-domain of $M$. Let $\chi$ be fixed a smooth cutoff function such that $\chi=1$ on $W$ and $\chi=0$ outside $W_1$.
We define
$$
\psi(x) =\exp_x(\chi(x)\cdot v(x)).
$$
Then $\psi$ is a diffeomorphism which equals to the identity map on the boundary and satisfies the following estimate
\be
\|\psi-Id\|_{C^{k-2}} \lesssim \|v\|_{C^{k-2}}\lesssim \|f_1\|_{C^{k-1}}.
\ee

This finishes the construction of $\psi$. We note that the diffeomorphism $\psi$ constructed above depends only on $\tilde{g}$ and the choice of the cutoff function $\chi$.

\medskip

Step 3.  Denote $\tilde{g}_2= \psi^*\tilde{g}_1$. It is clear that $\tilde{g}_2$ is isometric to $\tilde{g}_1$ and hence $\tilde{g}$. Moreover, $\tilde{g}_2$ has the same boundary normal coordinate as $g$ and hence $\tilde{g}$ also in $W$ (both $g$ and $\tilde{g}$ equal to $e$ there). Therefore, we can conclude that $\tilde{g}_2=\tilde{g}=e$ in $W$.

Let $\tilde{M}= M\backslash W$
and $f_2=\tilde{g}_2 -g$. Then the support of $f_2$ is contained in $\tilde{M}$.

We define $\varphi= \phi\circ\psi$. It is clear that $\tilde{g}_2= \psi^*\tilde{g}_1= \psi^* \phi^*\tilde{g}= \varphi^*\tilde{g}$.

\medskip

Step 4. We now have constructed two isometric copies of $\tilde{g}$ such that one is solenoidal with respect to the metric $g$ and the other is equal to $e$ in the neighborhood $W$ of $\p M$. Moreover, the following estimate holds for $f_1= \tilde{g}_1-g$ and $f_2= \tilde{g}_2-g$:
\begin{equation} \label{inequality1}
 \|f_2-f_1\|_{C^{l-3}} =  \|\tilde{g}_2-\tilde{g}_1\|_{C^{l-3}} \lesssim \|\psi -Id\|_{C^{l-2}}\lesssim
 \|f_1\|_{C^{l-1}}, \quad \mbox{for all}\,\, 3\leq l\leq k.
\end{equation}

By Proposition 1 in \cite{SU09}, we also have
\be \label{inequality2}
\|f_2^s -f_1\|_{C^{l-3}} \lesssim \|f_1\|_{C^{l-1}}^2 , \quad \quad \mbox{for all}\,\, 3\leq l\leq k.
\ee

\medskip

Step 5. Note that $\tilde{g}_2$ equals to $e$ in the neighborhood $W$ of $\p M$.
By Proposition \ref{prop1} and the fact that $\mathcal{H}^T(\tilde{g}_2)= \mathcal{H}^T(\tilde{g})$, we obtain
\be \label{a1}
\mathcal{H}^T(\tilde{g}_2)- \mathcal{H}^T(g) =\mathcal{H}^T(\tilde{g})- \mathcal{H}^T(g) \\
= \mathfrak{X}f_2 + O(\|f_2\|^2_{C^{5+[\f{d}{2}]}(M)})  \quad \mbox{in }  \,\, C^{3+[\f{d}{2}]}(S_{-}\p M).
\ee

Step 6. With the sub-manifold $\tilde{M}$ defined in Step 3, we let $\tilde{U}(x_j) \subset M$, $\alpha_j \in C_0^{\infty}(S_{-}\p  M)$, $j=1, 2... N$, and the
finite dimensional space $\mathfrak{L} \in \mathcal{S}L^2(S(\tau_2(M, \tilde{M}))$ be determined as in
Corollary \ref{thm-3}. By our assumption, $f_2 \perp \mathfrak{L}$, which further yields $f^s_2 \perp \mathfrak{L}$ by using the orthogonal decomposition of $f_2$. Therefore, the following estimate holds by Corollary \ref{thm-4}:
\be \label{a2}
 \|f^s_2\|_{H^{3+[\f{d}{2}]}(M)} \lesssim \sum_{j=1}^N \|\mathfrak{N}_{\alpha_j}f_2\|_{H^{2+[\f{d}{2}]}(\tilde{U}(x_j))}.
\ee

\medskip

Step 7. For each $j=1, 2, ...N$, we have
\bea
\|\mathfrak{N}_{\alpha_j}f_2\|_{H^{2+[\f{d}{2}]}(\tilde{U}(x_j))}&=&
\|\mathcal{L}^{\dag}\mathfrak{I}_{\alpha_j}^{\dag}\mathfrak{X}_{\alpha_j}f_2 \|_{H^{2+[\f{d}{2}]}(\tilde{U}(x_j))} \nonumber \\
&\lesssim & \|\mathfrak{I}_{ \alpha_j}^{\dag}\mathfrak{X}_{\alpha_j} f_2 \|_{H^{3+[\f{d}{2}]}(\tilde{U}(x_j))}\nonumber \\
&\lesssim & \|\mathfrak{I}_{\alpha_j}^{\dag}\mathfrak{X}_{ \alpha_j}  f_2 \|_{H^{3+[\f{d}{2}]}(M)} \nonumber \\
&\lesssim & \|\mathfrak{X}_{\alpha_j} f_2 \|_{H^{3+[\f{d}{2}]}(S_{-}\p M)} \label{a3}.
\eea

Step 8. By (\ref{a1}), (\ref{a2}) and (\ref{a3}), we deduce that
\beas
\|f^s_2\|_{H^{3+[\f{d}{2}]}(M)}
& \lesssim & \sum_{j=1}^N \|\mathfrak{X}_{{\alpha}_j} f_2 \|_{H^{3+[\f{d}{2}]}(S_{-}\p M)} \\
& \lesssim & \sum_{j=1}^N \|\alpha_j (\mathcal{H}^T(\tilde{g}) -\mathcal{H}^T(g) \|_{H^{3+[\f{d}{2}]}(S_{-}\p M)}
+ \|f_2\|_{C^{5+[\f{d}{2}]}(M)}^2.
\eeas
Therefore,
\beas
\|f_1\|_{C^{2}(M)} & \lesssim & \|f_2^s \|_{C^{2}(M)} + \|f_1\|_{C^{4}(M)}^2   \quad (\mbox{by taking $l=5$ in (\ref{inequality2})})\\
 & \lesssim & \|f_2^s \|_{H^{3+[\f{d}{2}]}(M)} + \|f_1\|_{C^{4}(M)}^2   \quad (\mbox{by Sobolev embedding theorem})\\
&\lesssim & \sum_{j=1}^N \|\alpha_j (\mathcal{H}^T(\tilde{g}) -\mathcal{H}^T(g) \|_{H^{3+[\f{d}{2}]}(S_{-}\p M)}
+ \|f_2\|_{C^{5+[\f{d}{2}]}(M)}^2 + \|f_1\|_{C^{4}(M)}^2  \\
& \lesssim & \sum_{j=1}^N \|\alpha_j (\mathcal{H}^T(\tilde{g}) -\mathcal{H}^T(g) \|_{H^{3+[\f{d}{2}]}(S_{-}\p M)}
+ \|f_1\|_{C^{7+[\f{d}{2}]}(M)}^2 \quad (\mbox{by taking $l=8+[\f{d}{2}]$ in (\ref{inequality1})}) \\
& \lesssim & \sum_{j=1}^N \|\alpha_j (\mathcal{H}^T(\tilde{g}) -\mathcal{H}^T(g) \|_{H^{3+[\f{d}{2}]}(S_{-}\p M)}
+ \|f_1\|_{C^{2}(M)} \cdot \|f_1\|_{C^{12+2[\f{d}{2}]}(M)}.
\eeas

Note that $\|f_1\|_{C^{12+2[\f{d}{2}]}(M)} \leq \|f_1\|_{C^{12+2[\f{d}{2}], \f{1}{2}}(M)} \lesssim \|f\|_{C^{13+2[\f{d}{2}], \f{1}{2}}(M)}$ (by taking $k=13+2[\f{d}{2}]$ in (\ref{inequality0}) for the last inequality).  By letting $\|f\|_{C^{13+2[\f{d}{2}], \f{1}{2}}(M)}$ sufficiently small, we obtain
\[
\|f_1\|_{C^{2}(M)}  \lesssim  \sum_{j=1}^N \|\alpha_j (\mathcal{H}^T(\tilde{g}) -\mathcal{H}^T(g) \|_{H^{3+[\f{d}{2}]}(S_{-}\p M)}.
\]
Finally, by taking $l=3$ in (\ref{inequality1}), we get the desired estimate
\[
\|f_2\|_{C^0(M)} \lesssim \|f_1\|_{C^{2}(M)}  \lesssim  \sum_{j=1}^N \|\alpha_j (\mathcal{H}^T(\tilde{g}) -\mathcal{H}^T(g) \|_{H^{3+[\f{d}{2}]}(S_{-}\p M)},
\]
which completes the proof of the theorem.

\section{Example of strong fold-regular metrics}  \label{sec-example}
In this section, we construct simple examples of strong fold-regular metrics. The main idea is the following: we perturb the Euclidean metric by creating a bump ``geometrically'' around some point. Consider the wave front 
generated by a source away from the bump. It is initially convex and spherical. After it passes through the bump, concavity can be developed. This concavity will eventually leads to conjugate points (or caustics). Since the metric is Euclidean except at the bump, explicit calculation can be carried out, and this enables us to find the explicit geometrical conditions which guarantee the strong fold-regular conditions. 


To begin with, we let $\varrho$ be a compactly supported smooth function of $\R^3$ whose support in contained in the unit ball. We consider the graph of the function $\varrho_z$: 
$$G_{\varrho, z}= \{(x, \varrho(x-z)); x \in \R^3\}.$$
$G_{\varrho, z}$ has a natural Riemannian metric which is induced from the Euclidean metric in $\R^4$, and which we denote by $g$.
Define
$
P: \R^3 \to G_{\varrho, z}
$
by 
$$
P (x) = (x, \varrho(x)).
$$
Then $P$ is a diffeomorphism  between $\R^3$ and $G_{\varrho, z}$. This diffeomorphism  
induces a metric $P^*g$ to $\R^3$, which is the one we want to construct. 

We now consider the Riemannian manifold $(\R^3, P^*g)$, which can be viewed as the flat space with one bump.  Let $x_0 $ be such that $|x_0 - z| \geq 1$ (here $|\cdot|$ means the Euclidean metric). We define $d_{P^*g}(x, x_0)$ to be the distance between the points $x$ and $x_0$ with respect to the metric $P^*g$. 

It is clear that level sets of the function $d_{P^*g}(\cdot, x_0)$
$$
\Gamma_{\beta} = \{x: d_{P^*g}(\cdot, x_0)= \beta \}
$$
are spheres for $\beta < |z-x_0|-1$. For $\beta > |z-x_0|-1$, only part of the level set is spherical. 
For properly chosen $\varrho$ and $z$, we assume that the following conditions are fulfilled: 

\begin{asump} \label{asump0-1}
For some $\beta_0 > |x-x_0|$, the map
$$
\exp_{x_0}: \{ v \in  T_{x_0} \R^3; \|v\|_{P^*g}= \beta_0\}  \to \Gamma_{\beta_0}
$$
is a diffeomorphism.  Moreover, the support of $\varrho(\cdot -z)$ (the bump)  is strictly contained in the bounded domain bounded by $\Gamma_{\beta_0}$. 
\end{asump}

For each $x\in \Gamma_{\beta_0}$, denote $N(x)$ the unit outward normal to the 
surface $\Gamma_{\beta_0}$.  Then $N=N(x)$ gives the Gaussian map of 
$\Gamma_{\beta_0}$. Let $(DN)_x$ be the differential of 
the Gaussian map $N$ at the point $x$. $(DN)_x$ is a linear and self-adjoint on the tangent space of $\Gamma_{\beta_0}$ at $x$. The two eigenvalues of this self-adjoint operator are the two principal curvatures of the surface  $\Gamma_{\beta_0}$ at $x$, which are denoted by $\kappa_1(x)$ and $\kappa_2(x)$ respectively. We also denote the corresponding eigenvectors by $U_1(x)$ and $U_2(x)$. 
We have 
$$
(DN)_x U_1(x)=  \kappa_1(x) U_1(x), \quad (DN)_x U_2(x)=  \kappa_1(x) U_2(x). 
$$

We further assume that 
 \begin{asump} \label{asump0-2}
The surface $\Gamma_{\beta_0}$ is concave at some point $x_1$. Moreover $$ \kappa_1(x_1) < \kappa_2(x_1) <0.$$
 \end{asump}
 

\begin{asump} \label{asump0-3}
The map 
$$
U_1: \Gamma_{\beta_0}  \to \cup_{x\in \Gamma_{\beta_0}} T_x \Gamma_{\beta_0} =\{ v\in \R^3; \|v\|_{P^*g}=1\}
$$
is a diffeomorphism in a neighborhood of $x_1$.  
\end{asump}

Now, we consider the diffeomorphism 
$
\exp_{x_0}: \{ \|v\|_{P^*g}= \beta_0; v \in T_{x_0} \R^3 \}  \to \Gamma_{\beta_0}
$. 
It induces a diffeomorphism $T: \{ v \in T_{x_0} \R^3 ; \|v\|_{P^*g}\geq   \beta_0 \}  \to \Gamma_{\beta_0} \times \R^{+}$:

$$
T(v) = ( \exp_{x_0} \beta_0 \f{v}{\|v\|_{P^*g}}, \|v\|_{P^*g} - \beta_0 ). 
$$
This diffeomorphism also yields a parameterization for the set $\{ v \in T_{x_0} \R^3 ; \|v\|_{P^*g} \geq   \beta_0 \}$. 

Define $F:\Gamma_{\beta_0} \times \R^{+} \to \R^3 $ by
$$
F(x, t) = x + t N(x). 
$$
Then we have
$$
\exp_{x_0} v = F \circ T (v),  \quad \|v\|_{P^*g} \geq   \beta_0. 
$$

It is clear that $v$ is a singular for $\exp_{x_0}$ if and only if the corresponding point $(x, t)=T(v)$ is singular for $F$. 
Note that $U_1(x)$, $U_2(x)$ and $N(x)$ forms an normal orthogonal basis for the space 
$T_{x}\R^3$. In this basis, we have 
$$
(DF)_{(x,t)}U_1(x) =  (1 +  t \kappa_1(x)) U_1(x) , \quad (DF)_{(x,t)}U_2(x) =  (1 +  t \kappa_2(x)) U_2(x), \quad (DF)_{(x,t)} N(x) = N(x). 
$$
As a result we see that 
$ \det  (DF)_{(x,t)} =0 $ if and only if 
$$
 1 +  t \kappa_1(x) = 0,  \quad \mbox{or}\,\,  1 +  t \kappa_2(x)=0. 
$$
Especially, we see that $(x_1, -\f{1}{\kappa_1(x_1)})$ is singular. In a sufficiently small neighbourhood of $(x_1, -\f{1}{\kappa_1(x_1)})$, the singular vectors of the map $F$ are given by
$
(x, -\f{1}{\kappa_1(x)}),
$
and the corresponding kernel space is the space spanned by the vector 
$
U_1(x). 
$

Therefore, we can conclude that Assumption \ref{asump0-2} implies that $T^{-1}(x_1, -\f{1}{\kappa_1(x_1)})$ is a fold type singular vector while Assumption \ref{asump0-3} implies that $T^{-1}(x_1, -\f{1}{\kappa_1(x_1)})$ is also strong-fold regular by Theorem 4.2 in 
\cite{SUAPDE}. 

The above construction can be quite general. In fact, for any given Riemannian metric,  let $\gamma$ be a geodesic of finite length without conjugate points. By perturbing the metric around some point on $\gamma$, we can create a conjugate point on the perturbed geodesic $\tilde{\gamma}$.  From the geometrical aspect of Assumption \ref{asump0-2}, we see that the created conjugate point is of fold type generically. It is strong fold-regular if Assumption \ref{asump0-3} is satisfied. We expect that Assumption \ref{asump0-3} is satisfied for generic surfaces. 

This completes the construction of strong-fold regular vectors. To construct a strong fold-regular point, we can apply the above construction consecutively to a properly selected set of geodesics emanated from $x_0$, we may construct a metric such that the two conditions in Definition \ref{def-fold-regular-point} are satisfied and thus $x_0$ is a strong fold-regular point.  Moreover, for this point $x_0$, one cannot find a ``complete'' set of geodesics which have no conjugate points. This shows that the strong fold-regular metrics introduced in this paper and \cite{BZ12}
are indeed non-trivial generalization of the ``regular'' metrics introduced in \cite{SU08AJM}.

\section*{Appendix A: Relations between the geodesic flow and the lens relation}

We present the relations between the geodesic flow and the lens relation in the form (\ref{s-relation}).
Let $\mathcal{H}_0^t$ denote the geodesic flow induced by the Euclidean metric in $\R^d$ and let $g$ be
the background metric. Then the following identity holds:
\be
\mathcal{H}^T(g)(x_0, \xi_0)= \mathcal{H}_0^{T- L(g)(x_0, \xi_0)}(\Sigma(g)(x_0, \xi_0)), \quad
(x_0, \xi_0)\in S_{-}\p M.  \label{identity3-1-1}
\ee

On the other hand, for each $(x, \xi)\in S(\R^d\backslash M)$, we define $\kappa(x, \xi)$ to be the first positive moment when the orbit of the geodesic flow $\mathcal{H}_0^t(x, -\xi)$ hits the boundary $\p M$. By the assumptions on the background metric $g$, $\kappa$ is well-defined and is smooth in a sufficiently small neighborhood of $\mathcal{H}^T(g)(x_0, \xi_0)$ for each
$(x_0, \xi_0)\in S_{-}\p M$. As a result, the following identities hold for all $(x_0, \xi_0)\in S_{-}\p M$:
\bea
L(g)(x_0, \xi_0)&=& \kappa (\mathcal{H}^T(g)(x_0, \xi_0)); \label{identity3-1-2} \\
\Sigma(g)(x_0, \xi_0)& =& \mathcal{H}_0^{ \kappa (\mathcal{H}^T(g)(x_0, \xi_0))}(\mathcal{H}^T(g)(x_0, \xi_0)). \label{identity3-1-3}
\eea

Now, by taking derivatives with respect to $g$ in the above identities (\ref{identity3-1-1}),
(\ref{identity3-1-2}) and (\ref{identity3-1-3}), we obtain the following result.

\begin{lem} \label{lem-relation2}
There exist two smooth matrix-valued functions $B \in C^{\infty}(S_{-}\p M)$ and $C \in C^{\infty}(S_{-}\p M)$ such that for all $(x_0, \xi_0)\in S_{-}\p M$,
\[
\f{\delta \mathcal{H}^T(g)}{\delta \,g}(x_0, \xi_0) = B(x_0, \xi_0)\left(\begin{array}{l}
 \f{\delta \Sigma(g)}{\delta \,g}(x_0, \xi_0) \\
  \f{\delta L(g)}{\delta \, g}(x_0, \xi_0) \end{array}\right);
 \quad \left(\begin{array}{l}
 \f{\delta \Sigma(g)}{\delta g}(x_0, \xi_0) \\
  \f{\delta L(g)}{\delta g}(x_0, \xi_0)\end{array}\right) = C(x_0, \xi_0) \f{\delta \mathcal{H}^T(g)}{\delta \, g}(x_0, \xi_0).
\]
\end{lem}

\begin{rmk} \label{rmk-1}
 The matrix-valued function $C$ is only smooth in $S_{-}\p M$ and it may not be extended smoothly to $\overline{S_{-}\p M}$. This is because the function $\kappa$ may not have a smooth extension to
 $\mathcal{H}^T(g)(\overline{S_{-}\p M})$.
\end{rmk}

\section*{Appendix B: Proof of Proposition \ref{prop1}}

We first linearize the nonlinear operator which maps
a metric to its corresponding Christoffel symbol at $g$.
We denote by $\Gamma(g+f)$ the Christoffel symbol corresponding to the metric $g+f$.
Observe that
\be \label{christoffel}
\Gamma_{ij}^k(g)= \f{1}{2}\left( \f{\p g_{jp}}{\p x^i} + \f{\p g_{ip}}{\p x^j}-\f{\p g_{ij}}{\p x^p} \right).
\ee
We can show that
$$
\mathcal{L}f=\f{\delta \, \Gamma(g)}{\delta \, g} f.
$$

We now present an estimate for the quality of the linear approximation.
\begin{lem} \label{lem-linear1}
For $f\in S(\tau_2M)$ with sufficiently small $C^2$ norm, the following estimate holds for any positive integer $k$.
\[
 \|\Gamma(g+f) - \Gamma(g) - \mathcal{L} f \|_{C^{k}(M)} \lesssim \|f\|^2_{C^{k+1}(M)}.
\]
\end{lem}

{\bf{Proof}}. For fixed $g$ and $f$, define $\Gamma(\tau)= \Gamma(g+ \tau f)$. It is clear that
 \[
 \Gamma(g+f)- \Gamma(g) = \Gamma(1)- \Gamma(0) = \Gamma'(0) + \f{1}{2}\Gamma''(\eta), \quad \mbox{for some} \,\, \eta \in [0, 1]
 \]
and $\Gamma'(0)= \mathcal{L}f$.
To prove the lemma, it suffices to show that
 $$
  \|\Gamma''(\tau)\|_{C^{k}} \lesssim \|f\|_{C^{k+1}}^2.
 $$
We argue as follows. First, from Formula (\ref{christoffel}), we can deduce that
\[
2\, \Gamma_{ij}^k(\tau) g_{kl}(\tau) = \f{\p g_{jl}}{\p x^i}+\f{\p g_{il}}{\p x^j}-\f{\p g_{ij}}{\p x^l} + \tau (\f{\p f_{jl}}{\p x^i}+\f{\p f_{il}}{\p x^j}-\f{\p f_{ij}}{\p x^l}).
\]
Next, we differentiate the above identity with respect to $\tau$ to obtain
\be \label{identity111}
 2\, \dot{\Gamma}_{ij}^k(\tau) g_{kl}(\tau) + 2 \,\Gamma_{ij}^k(\tau) f_{kl} = \f{\p f_{jl}}{\p x^i}+\f{\p f_{il}}{\p x^j}-\f{\p f_{ij}}{\p x^l}.
\ee
We then multiply both sides of the above identity by $g^{lp}$ and summate over $l$. This gives
$$
2\, \dot{\Gamma}_{ij}^k(\tau) = - 2 \,\Gamma_{ij}^k(\tau) f_{kl} g^{lp} +  g^{lp}(\f{\p f_{jl}}{\p x^i}+\f{\p f_{il}}{\p x^j}-\f{\p f_{ij}}{\p x^l}).
$$
It follows that
$$
\|\dot{\Gamma}_{ij}^k(\tau)\|_{C^k} \lesssim \|f\|_{C^{k+1}}.
$$
Now, differentiating the identity (\ref{identity111}) with respect to $\tau$ again,
$$
 \ddot{\Gamma}_{ij}^k(\tau) g_{kl}(\tau) + 2\,\dot{\Gamma}_{ij}^k(\tau) f_{kl}=0.
$$
Using the same trick as we did for $\dot{\Gamma}_{ij}^k(\tau)$, we can conclude that
$$
\|\ddot{\Gamma}_{ij}^k(\tau)\|_{C^k} \lesssim \|\dot{\Gamma}_{ij}^k(\tau)\|_{C^k} \|f\|_{C^k} \lesssim \|f\|_{C^{k+1}}^2,
$$
which yields the desired estimate and completes the proof of the lemma.

\medskip

We next linearize the nonlinear operator which maps
a Christoffel symbol to its corresponding geodesic flow at $\Gamma$. By a similar
argument as in \cite{BZ12}, we can show that the following holds.

\begin{lem} \label{lem-linear2}
Let $k$ be a nonnegative integer, the following estimate holds for the linearization of the nonlinear operator which maps
Christoffel symbols to its corresponding geodesic flow at $\Gamma$:
\be
  \| \mathcal{H}^T(\Gamma + \tilde{\Gamma}) - \mathcal{H}^T(\Gamma) -\mathfrak{I} (\tilde{\Gamma}) \|_{C^k(S_{-}\p M)}
  \lesssim \| \tilde{\Gamma}\|^2_{C^{k+1}(M)}.
\ee
\end{lem}

{\bf{Proof of Proposition \ref{prop1}}}: It is a consequence of Lemma \ref{lem-linear1}, \ref{lem-linear2} and
the chain rule.

\bigskip
\noindent{\bf{Conflict of Interest}}: The authors declare that they have no conflict of interest.

\end{document}